\documentclass[final,10pt]{amsart}

\usepackage[left=1.3in, right=1.3in, top=1.5in, bottom=1.5in]{geometry}

\usepackage[foot]{amsaddr} 

\usepackage[]{showkeys}

\usepackage{verbatim}
\usepackage{amssymb}
\usepackage{amsthm}
\usepackage{amsmath}
\usepackage{amsfonts}
\usepackage{latexsym}
\usepackage{mathtools} 
\usepackage{enumitem} 
\usepackage[english]{babel} 
\usepackage{mathabx}

\usepackage{float}

\usepackage[usenames,dvipsnames]{xcolor} 
\definecolor{dkblue}{RGB}{30,90,140} 
\usepackage[colorlinks=true, pdfstartview=FitV, linkcolor=dkblue, citecolor=dkblue, urlcolor=dkblue]{hyperref}    

\definecolor{mydarkbluett}{RGB}{12,111,174}

\usepackage{todonotes}

\emergencystretch=\maxdimen
\theoremstyle{plain}
\newtheorem{thm}{Theorem}
\newtheorem{remark}[thm]{Remark}

\numberwithin{equation}{section}
\numberwithin{thm}{section}

\newcommand{\pv}{\text{pv}}
\newcommand{\pa}{\partial_{\alpha}}

\newcommand{\minspace}{\hspace{0.05cm}}

\newcommand{\R}{\mathbb R}

\newcommand{\bm}[1]{\boldsymbol{#1}}

\DeclareMathOperator{\Lip}{Lip}

\begin{document}
	

	\title[Quantitative H\"older Estimates for Even S.I.O. on patches]{Quantitative H\"older Estimates  for \\ Even Singular Integral Operators on patches}
	
	\author[F. Gancedo]{Francisco Gancedo$^\dagger$}
	\address{$^\dagger$Departamento de An\'{a}lisis Matem\'{a}tico $\&$ IMUS, Universidad de Sevilla, C/ Tarfia s/n, Campus Reina Mercedes, 41012 Sevilla, Spain.  \href{mailto:fgancedo@us.es}{fgancedo@us.es}}
	
	\author[E. Garc\'ia-Ju\'arez]{Eduardo Garc\'ia-Ju\'arez$^{\ddagger}$}
	\address{$^{\ddagger}$Departament de Matemàtiques i Informàtica, Universitat de Barcelona, Gran Via de les Corts Catalanes, 585 08007, Barcelona, Spain. \href{mailto:egarciajuarez@ub.edu}{egarciajuarez@ub.edu}}

	
	\begin{abstract}
		In this paper we show a constructive method to obtain $\dot{C}^\sigma$ estimates of even singular integral operators on characteristic functions of domains with $C^{1+\sigma}$ regularity, $0<\sigma<1$. 
		This kind of functions were shown in first place to be bounded (classically only in the $BMO$ space) to obtain global regularity for the vortex patch problem \cite{Chemin1993, BertozziConstantin1993}. This property has then been applied to solve different type of problems in harmonic analysis and PDEs. Going beyond in regularity, the functions are discontinuous on the boundary of the domains, but $\dot{C}^{\sigma}$ in each side.   
		This $\dot{C}^{\sigma}$ regularity has been bounded by the $C^{1+\sigma}$ norm of the domain \cite{Depauw1999, MateuOrobitgVerdera2009,Mitrea2Verdera2016}. Here we provide a quantitative bound linear in terms of the $C^{1+\sigma}$ regularity of the domain.  This  estimate shows explicitly the dependence of the lower order norm and the non-self-intersecting property of the boundary of the domain. As an application, this  quantitative estimate is used in a crucial manner to the free boundary incompressible Navier-Stokes equations providing new global-in-time regularity results in the case of viscosity contrast \cite{GancedoG-J2021}.    
		
	\end{abstract}
	
	\setcounter{tocdepth}{2}
	
	\maketitle
	\tableofcontents

	\section{Introduction}
	
	In this paper we deal with characteristic functions of domains
	$$
	1_{ D}(x)=\left\{\begin{array}{ll}
		1,& x\in  D,\\
		0,& x\in \R^n\smallsetminus\overline{ D}.
	\end{array}\right.
	$$
	The main interest is to study singular integral operators of Calderón-Zygmund type applied on this kind of functions and obtained as follows
	\begin{equation}
		\label{T1O}
		T(1_ D)(x)=\pv\int_{\R^n}K(x-y)1_{ D}(y) dy=\lim_{\varepsilon\to 0^+}\int_{|x-y|>\varepsilon}K(x-y)1_{ D}(y) dy,\qquad x\in \R^n.  
	\end{equation}
	Above, $\pv$ stands for principal value and the kernel $K$ is homogeneous of degree $-n$, given by
	$$
	K(x)=\frac{\Omega(x)}{|x|^n},\quad\Omega(\lambda x)=\Omega(x)\quad \forall \lambda>0,\quad \int_{|x|=1}\Omega(x)d\sigma(x)=0.
	$$
	In the classical theory of singular integrals, the function $T(1_ D)(x)$ belongs to the $BMO$ space \cite{Stein93}. In particular, for odd kernels, it is not difficult to show that as $x$ approaches to a point on $\partial D$ the function $T(1_ D)(x)$ is not bounded and diverges to infinity logarithmically. On the other hand, when the kernel is even,
	$$
	\Omega(x)=\Omega(-x), 
	$$
	a new geometric cancellation was found in \cite{Chemin1993,BertozziConstantin1993} which  shows that $T(1_ D)(x)$ belongs to $L^\infty$. This $L^\infty$ bound is given in terms of the $C^{1+\sigma}$ norm of the domain, $0<\sigma<1$. The regularity of the boundary together with the fact that the kernel has mean zero on half spheres cancel the singularity on the boundary of the domain. The motivation was to show preservation of $C^{1+\sigma}$ regularity for domains moving by the 2D Euler equations; i.e. global in time existence for the vortex patch problem \cite{Chemin1993,BertozziConstantin1993}. 
	
	From the harmonic analysis point of view, Calder\'on-Zygmund operators with smooth and even kernel have been studied specifically as they satisfy stronger inequalities than general ones. In  \cite{MOV2011}, it is shown that the following pointwise inequality holds for even, higher-order Riesz Transforms 
	$$
	T^*f(x)\leq CM(Tf)(x),
	$$
	where $T^*$ is the maximal singular integral and $M$ is the Hardy-Littlewood maximal operator. It yields a stronger estimate than the classical Cotlar's inequality \cite{Torchinsky1986}.

	The extra cancellation providing $L^\infty$ bounds has been extensively used in different PDEs problems. Considering the Beltrami equation, it guarantees that the solutions are bi-Lipschitz  \cite{MateuOrobitgVerdera2009}. For the Muskat problem, modeling the evolution of incompressible immiscible fluids in porous media or Hele-Shaw cells, this bound yields lack of squirt singularities \cite{CordobaGancedo2010} (also known in the literature as splat singularities). For multidensional aggregation equations with a Newtonian potential, it provides propagation of $C^{1+\sigma}$ regularity up to the blow-up \cite{BGLV2016}. In the two dimensional inhomogeneous Navier-Stokes equations modeling the evolution of incompressible fluids of different densities, this $L^\infty$ bound provides global-in-time regularity for higher order norms ($W^{2,\infty}$ and $C^{2+\sigma}$) of the moving free boundary between the fluids \cite{GancedoG-J2018}. In \cite{GancedoG-J2020}, a combination of parabolic and elliptic estimates together with this $L^\infty$ bound are used to propagate the same higher order norms but for Boussinesq temperature fronts. See also \cite{CMOV2021} for recent developments in contour dynamics for non-linear transport equations. In all these results the singular integral operators are given with $\Omega$ a polynomial function.
	
	In this work, we go further in order to control higher regularity for functions given by \eqref{T1O} with  even kernel. Despite the fact that these functions are discontinuous on $\partial D$, it is possible to obtain $C^{\sigma}$ regularity in $ D$ and in $\R^n\smallsetminus\overline{ D}$. In \cite{MateuOrobitgVerdera2009, Mitrea2Verdera2016}, this regularity has been shown together with qualitative bounds of the form
	\begin{equation}\label{VerderaBound}
		\|T(1_ D)\|_{C^\sigma(\overline{ D})}\leq C P (\| D\|_{C^{1+\sigma}}),
	\end{equation}
with $ P $ a polynomial and $C>0$ a constant depending on $n$, $\sigma$, and the geometry of the domain $D$. Above, the function $T(1_ D)$ is extended continuously on $\overline{ D}$. The latter paper also characterizes the regularity of the domains in terms of odd singular integrals operators on $\partial D$. It uses harmonic analysis techniques and Clifford algebras as a generalization of the field of complex numbers to higher dimensions.
	In this paper, we show that the bound above can be improved to make it linear in the higher regularity norm of the boundary. Moreover, the dependency on the arc-chord condition is made explicit:
	\begin{equation}\label{sharpc1g}
		    \begin{aligned}
		        \|T(1_D)\|_{\dot{C}^{\sigma}(\overline{ D})\cup \dot{C}^{\sigma}(\R^n\smallsetminus D)}\leq 
			C(1\!+\!|\partial D|)  P (\|D\|_{*}\!+\!\|D\|_{\Lip})(1+\|D\|_{\dot{C}^{1+\sigma}}).
		    \end{aligned}
		\end{equation}
	Above, $C=C(n,\sigma)$, $\Lip$ stands for Lipschitz, $\| D\|_*$ measures the non self-intersecting property of the boundary $\partial D$, $ P $ is a polynomial function, and $|\partial D|$ denotes the $(n-1)$-dimensional surface area. The higher order norm is homogeneous, given for a function by
	$$
	\|f\|_{\dot{C}^{1+\sigma}}=\|\nabla f\|_{\dot{C}^{\sigma}}=\sup_{x\neq y}\frac{|\nabla f(x)-\nabla f(y)|}{|x-y|^{\sigma}},\quad 0<\sigma<1.
	$$
	See below for more details about the notation.
	During the review process of this article the referee pointed out work \cite{Depauw1999}, in which the author also proves an estimate similar to the one above. However, our proof is different, working at the level of the interface via contour dynamics methods.
	
	An important motivation of these estimates comes from a classical two dimensional fluid mechanics problem. Concretely, the dynamics of two incompressible immiscible fluids evolving by the inhomogeneous Navier-Stokes equations. In that problem, the viscosity can be understood as a patch function and the gradient of the velocity is related to the viscosity by combinations of second and fourth-order Riesz transforms,
	$$
	T=\partial_j\partial_k(-\Delta)^{-1}(I-\nabla(-\Delta)\nabla\cdot),\quad j,k=1,...,n,\quad I \mbox{ the identity.}
	$$
    In \cite{GancedoG-J2021}, global-in-time well-posedness for the evolution of $C^{1+\sigma}$  interfaces between the two fluids is proved. Global-in-time regularity was recently shown in \cite{PaicuZhang2020} for $H^{5/2}$ Sobolev regularity of the interface instead of $C^{1+\sigma}$ and by using striated regularity.	In the argument of the proof in \cite{GancedoG-J2021}, the estimate \eqref{sharpc1g} is used in an important manner. In particular,  we emphasize the importance of the quantitative bound of the non self-intersection condition  $\|D\|_{*}$. This quantity has to be controlled globally, since it is known that free-boundary incompressible Navier-Stokes can developed finite-time pointwise particle collision on the free interface \cite{CCFGG19,CS19}. 
	

	\vspace{0.5cm}

The rest of the paper is structured as follows.  Section \ref{main_result} contains the statement of the main results: Theorems \ref{MainTheorem} and \ref{cor1}. It describes how the operators \eqref{T1O} can be studied in terms of odd operators on the boundary, yielding Theorem \ref{cor1} as a corollary of Theorem \ref{MainTheorem}. The rest of the paper, Section \ref{main_thm_sec}, is the proof of Theorem \ref{MainTheorem}. To study the H\"older regularity of the operators involved, the proof distinguishes three situations: when the two points are on the boundary (Subsection \ref{on_boundary}), \textit{near} the boundary (Subsection \ref{near_boundary}), and \textit{far} from the boundary (Subsection \ref{reg_far}). The deciding cut-off is defined in terms of the non self-intersecting condition and the $\dot{C}^{1+\sigma}$ regularity of the domain. In the second scenario, we need to consider further whether the separation between the points occurs \textit{mostly} in normal or tangential direction.  The \textit{nearly} normal direction case is decomposed in purely normal (Case 1) and tangential (Case 2) differences, each one estimated through delicate splittings of the singular integrals. The \textit{nearly} tangential case is reduced, using a fixed point argument, to the purely normal plus on the boundary cases. Finally, the third situation is less singular.

	\section{Main Result}\label{main_result}
	
	We consider higher-order Riesz transform operators of even order $2l$, $l\geq1$. That is, we deal with Calder\'on-Zygmund operators given by
	\begin{equation}\label{calderonzygmund}
		\begin{aligned}
			T(f)(x)= \lim_{\varepsilon\to 0}\int_{|x-y|>\varepsilon} K(x-y)f(y)dy,
		\end{aligned}
	\end{equation}
	where
	\begin{equation}\label{kernel}
		K(x)=\frac{P_{2l}(x)}{|x|^{n+2l}},
	\end{equation}
	and $P_{2l}(x)$ is a homogeneous harmonic polynomial of degree $2l$ in $\R^n$.
	We want to study the H\"older regularity of the operator $T$ applied to the characteristic function $1_ D(x)$ of a $C^{1+\sigma}$ domain $ D$,
	\begin{equation}\label{TOmega}
		T(1_ D)(x)=\pv\int_{\R^2}K(x-y)1_{ D}(y) dy=\pv\int_{ D}K(x-y) dy,\qquad x\in \R^n.
	\end{equation}
	We recall that the operators \eqref{calderonzygmund} have explicit Fourier multipliers, 
	\begin{equation}\label{multip}
		\mathcal{F}\Big(\frac{P_{m}(x)}{|x|^{n+m-\alpha}}\Big)(\xi)=c_{m,\alpha,n}\frac{P_{m}(\xi)}{|\xi|^{m+\alpha}},
	\end{equation}
	with $0<\alpha\leq n$ and
	\begin{equation*}
		c_{m,\alpha,n}=i^m\pi^{\frac{n}2-\alpha}\frac{\Gamma\big(\frac{m+\alpha}{2}\big)}{\Gamma\big(\frac{m+n-\alpha}{2}\big)}.
	\end{equation*}
	\begin{remark}\label{harmonic}
		Any homogeneous polynomial $P_k$ of degree $k$ can be written as $P_k(x)=p_{k}(x)+|x|^2p_{k-2}(x)$, where $p_k$ is a homogeneous harmonic polynomial of degree $k$ and $p_{k-2}$ is homogeneous of degree $k-2$ (Sec. 3 in Chapter 3  of \cite{Stein1970}). Thus the restriction to harmonic polynomials in \eqref{kernel} involves no loss of generality.
	\end{remark}
	Using Euler's homogeneous function theorem and integration by parts, the regularity of \eqref{TOmega} can be studied through the associated operators 
	\begin{equation}\label{Soperator}
		S(f)(x)=\pv \int_{\partial D} k(x-y) f(y)dS(y),\qquad x\in\mathbb{R}^n,
	\end{equation}
	where the kernel $k(x)$ is given by
	\begin{equation*}
		k(x)=\frac{Q_{2l-1}(x)}{|x|^{n+2l-2}},
	\end{equation*}
	and $Q_{2l-1}$ is a homogeneous harmonic polynomial of degree $2l-1$.
	
Since one of the main motivations of these results are physical, we show the techniques in dimension three. An analogous approach provides the proof in any dimension. 
	
	We will denote by $D$ a non-self-intersecting bounded domain of class $C^{1+\sigma}$. It is defined as follows. Denote by $V_j$, $j = 1,...,J$, the neighborhoods that provide local charts of the boundary $\partial D$ in such a way that for any $x\in\partial D$ there exists a $V_j\subset\R^2$ such that $x=Z(\alpha)$ with $\alpha\in V_j$, with well-defined normal vector. 
	To measure the non-self intersection and the regularity of the parameterization, we define
	$$
	\|D\|_*=\|F(Z)\|_{L^\infty}=|\partial_\alpha Z|_{\inf}^{-1}<\infty,
	$$
	where
	$$
	|\pa Z|_{\inf}=\min_{j=1,...,J}\{\min \{\inf_{\alpha\in V_j}|\partial_{\alpha_1}Z(\alpha)|,\inf_{\alpha\in V_j}|\partial_{\alpha_2}Z(\alpha)|,\inf_{\alpha\neq\beta, \alpha,\beta\in V_j}\frac{|Z(\alpha)-Z(\beta)|}{|\alpha-\beta|}\}\}.
	$$
	The Lipschitz norm is given by
	$$
	\|D\|_{\Lip}=\max_{j=1,...,J}\sup_{\alpha\neq\beta, \alpha,\beta\in V_j}\frac{|Z(\alpha)-Z(\beta)|}{|\alpha-\beta|}<\infty,
	$$
	and the H\"older seminorm by 
	$$
	\|D\|_{\dot{C}^{1+\sigma}}=\max_{j=1,...,J}\sup_{\alpha\neq\beta, \alpha,\beta\in V_j}\frac{|\nabla Z(\alpha)-\nabla Z(\beta)|}{|\alpha-\beta|^{\sigma}}<\infty.
	$$
		There exists also a well-defined normal vector given by $$N=\partial_{\alpha_1}Z\wedge \partial_{\alpha_2}Z.$$
	For convenience, we take the parametrization so that $N$ is pointing towards the interior of the surface. 
 Now we are in position to state our main theorem:

	\begin{thm}\label{MainTheorem}
		Assume $D$ is a bounded domain of class $C^{1+\sigma}$, $0<\sigma<1$. Then, the operator \eqref{Soperator} maps boundedly $\dot{C}^{\sigma}(\partial D)$ into $\dot{C}^{\sigma}(\overline{ D})\cup \dot{C}^{\sigma}(\R^n\smallsetminus D)$. Moreover, the following bound holds
		\begin{equation*}
			\|S(f)\|_{\dot{C}^{\sigma}(\overline{ D})\cup \dot{C}^{\sigma}(\R^n\smallsetminus D)}\leq C(1\!+\!|\partial D|)P(\|D\|_{*}\!+\!\|D\|_{\Lip})\Big(\|f\|_{C^{\sigma}}+\|f\|_{L^\infty}\|D\|_{\dot{C}^{1+\sigma}}\Big),
	\end{equation*}
with $P$ a polynomial function depending on $S$, and $C=C(n,\sigma)$.
	\end{thm}

	As indicated before, we can write \eqref{TOmega} as a sum of terms of the form \eqref{Soperator} with $f=N_j$. Therefore, we have the following result:
	
	\begin{thm}\label{cor1}
		Assume $D$ is a bounded domain of class $C^{1+\sigma}$, $0<\sigma<1$. Then, the Calder\'on-Zygmund operator \eqref{calderonzygmund} applied to the characteristic function of $ D$, \eqref{TOmega}, defines a piecewise $\dot{C}^{\sigma}$ function, \begin{equation*}
			T(1_D)\in \dot{C}^{\sigma}(\overline{D})\cup \dot{C}^{\sigma}(\R^n\smallsetminus D).
		\end{equation*}
		Moreover, it satisfies the bound
		\begin{equation*}
		    \begin{aligned}
		        \|T(1_D)\|_{\dot{C}^{\sigma}(\overline{ D})\cup \dot{C}^{\sigma}(\R^n\smallsetminus D)}\leq 
			C(1\!+\!|\partial D|) P(\|D\|_{*}\!+\!\|D\|_{\Lip})(1+\|D\|_{\dot{C}^{1+\sigma}}),
		    \end{aligned}
		\end{equation*}
with $P$ a polynomial function depending on $T$, and $C=C(n,\sigma)$.
	\end{thm}

	

	\section{Proof of Theorem \ref{MainTheorem}}\label{main_thm_sec}

	Without loss of generality, we show the proof for the following case 
	\begin{equation*}
		\begin{aligned}
			S(f)(x)=\pv\int_{\partial D} k(x-y)f(y)dS(y),      \end{aligned}
	\end{equation*}
	with 
	\begin{equation*}
		k(x)=\frac{x_1x_2x_3}{|x|^5}.
	\end{equation*}
	The developed techniques can be applied to any other odd homogeneous polynomial and any dimension. We choose this case to show more clearly the crucial steps. The case of $k(x)$ of degree one is more direct. If the degree is greater than three the approach is the same but technically longer. With the kernel chosen, we  provide a constructive and direct method, showing the main difficulties and cancellations.

	\subsection{Regularity on the boundary.}\label{on_boundary} First, consider an atlas of the surface $\partial D$ and, on a given chart, fix a cut-off $\eta>0$ and define the ball
	$$
	A_\eta=\{y\in\partial D:|x-y|<\eta\}.
	$$
	Consider any $h\in \R^3$ such that
	\begin{equation*}
		|h|\leq \overline{\eta}=\frac{\eta}{4(1+\|F(Z)\|_{L^\infty})(1+\|\pa Z\|_{L^\infty})},
	\end{equation*}
	then $x+h\in A_{\overline{\eta}}$, and we will generally write 
	\begin{equation*}
		x=Z(\alpha),\qquad x+h=Z(\beta).    
	\end{equation*}
	We will also use the middle coordinate point
	\begin{equation}\label{mid_alphabeta}
	    \xi=\frac{\alpha+\beta}{2}.
	\end{equation}
	Then, 
	\begin{equation}\label{holdersplit_boundary}
		\begin{aligned}
			S(f)(x)-S(f)(x+h)&=\pv\int_{A_\eta}\Big(k(x-y)-k(x+h-y)\Big)f(y)dS(y)\\
			&\quad+\int_{\partial D \smallsetminus A_\eta}\Big(k(x-y)-k(x+h-y)\Big)f(y)dS(y)\\
			&=I+II.
		\end{aligned}
	\end{equation}
	The second term is away of the singular part and thus more regular, 
	\begin{equation*}
		\begin{aligned}
			|II|\leq \frac{C}{\eta^3}|\partial D|\|f\|_{L^\infty}|h|. 
		\end{aligned}
	\end{equation*}
	The first term is given by
	\begin{equation*}
		\begin{aligned}
			I&=\int_{Z^{-1}(A_\eta)} \big(k(Z(\alpha)-Z(\gamma))-k(Z(\beta)-Z(\gamma))\big)f(Z(\gamma))|N(\gamma)|d\gamma.
		\end{aligned}    
	\end{equation*}
	For simplicity of notation, we will denote
	\begin{equation*}
		g(\gamma)=f(Z(\gamma))|N(\gamma)|,
	\end{equation*}
	so we have that
	\begin{equation}\label{gbound}
	    \begin{aligned}
	        \|g\|_{L^\infty}&\leq \|\pa Z\|_{L^\infty}^2\|f\|_{L^\infty},\\
	        \|g\|_{\dot{C}^{\sigma}}&\leq \|f\|_{\dot{C}^\sigma}\|\pa Z\|_{L^\infty}^{2+\sigma}+\|\pa Z\|_{\dot{C}^\sigma}\|\pa Z\|_{L^\infty}\|f\|_{L^\infty}.
	    \end{aligned}
	\end{equation}
	The nonlinear kernels in $I$ are neither odd nor given by a derivative. We decompose the $I$ term as follows
	\begin{equation}\label{Isplit_boundary}
		\begin{aligned}
			I=I_1+I_2,
		\end{aligned}
	\end{equation}
	where
    \begin{equation*}
		\begin{aligned}
			I_1=&\int_{Z^{-1}(A_\eta)}\Big(k(Z(\alpha)-Z(\gamma))-k(\pa Z(\xi)(\alpha-\gamma))\Big)g(\gamma)\minspace d\gamma\\
			&-\int_{Z^{-1}(A_\eta)}\Big(k(Z(\beta)-Z(\gamma))-k(\pa Z(\xi)(\beta-\gamma))\Big)g(\gamma)\minspace d\gamma,   
		\end{aligned}
	\end{equation*}
	and
	\begin{equation*}
		\begin{aligned}
			I_2&=\int_{Z^{-1}(A_\eta)}\Big(k(\pa Z(\xi)(\alpha-\gamma))-k(\pa Z(\xi)(\beta-\gamma))\Big)g(\gamma)\minspace d\gamma.
		\end{aligned}
	\end{equation*}
	We estimate $I_2$ first. Since the kernel is odd, we first isolate the singularity from the boundary.
		We notice that 
	\begin{equation}\label{distances}
		|\alpha-\beta|\leq \|F(Z)\|_{L^\infty}|h|\leq \|F(Z)\|_{L^\infty}\overline{\eta}\leq \frac{\eta}{4\|\pa Z\|_{L^\infty}},
	\end{equation}
	and
	\begin{equation*}
		d(\alpha, \partial Z^{-1}(A_{\eta}))\geq\frac{\eta}{\|\pa Z\|_{L^\infty}},\quad d(\beta, \partial Z^{-1}(A_\eta))\geq\frac{3\eta}{4\|\pa Z\|_{L^\infty}},
	\end{equation*}
	thus, we take a smooth cut-off $\chi(\gamma)$ defined as $\chi(\gamma)=1$ for $|\alpha-\gamma|\leq \frac{\eta}{2\|\pa Z\|_{L^\infty}}$, $\chi(\gamma)=0$ for $|\alpha-\gamma|\geq\frac{3\eta}{4\|\pa Z\|_{L^\infty}}$, and radial centered at $\alpha$. 
	Introducing the cut-off in $I_{2}$, we have
	\begin{equation}\label{I21split_boundary}
		\begin{aligned}
			I_{2}&=I_{2,1}+I_{2,2},
		\end{aligned}
	\end{equation}
	with
	\begin{equation*}
		\begin{aligned}
			I_{2,1}&=\int_{Z^{-1}(A_\eta)}\Big(k(\pa Z(\xi)(\alpha-\gamma))-k(\pa Z(\xi)(\beta-\gamma))\Big)\chi(\gamma)g(\gamma)\minspace d\gamma,\\
			I_{2,2}&=\int_{Z^{-1}(A_\eta)}\Big(k(\pa Z(\xi)(\alpha-\gamma))-k(\pa Z(\xi)(\beta-\gamma))\Big)(1-\chi(\gamma))g(\gamma)\minspace d\gamma.
		\end{aligned}
	\end{equation*}
	Along the paper, we will need to control the singularity in the kernel $k(\pa Z(\xi)(\alpha-\gamma))$. We will use that it is comparable to $|\alpha-\gamma|^{-2}$. In fact, since the domain is regular and of class $C^{1+\gamma}$, it holds that
	\begin{equation}\label{aux}
	   |\partial_{\alpha_1}Z(\xi)\cdot\partial_{\alpha_2}Z(\xi)|\leq (1-\varepsilon)|\partial_{\alpha_1}Z(\xi)||\partial_{\alpha_2}Z(\xi)|,
	\end{equation}
	for some $\varepsilon>0$. Therefore, 
	\begin{equation*}
	    \begin{aligned}
	        |\partial_\alpha Z(\xi)(\alpha-\gamma)|^2&\geq (\alpha_1-\gamma_1)^2|\partial_1 Z(\xi)|^2+(\alpha_2-\gamma_2)^2|\partial_2 Z(\xi)|^2\\
	        &\quad-2(1-\varepsilon)|\alpha_1-\gamma_1||\alpha_2-\gamma_2||\partial_1 Z(\xi)||\partial_2 Z(\xi)|\\
	        &\geq \varepsilon\Big((\alpha_1-\gamma_1)^2|\partial_1 Z(\xi)|^2+(\alpha_2-\gamma_2)^2|\partial_2 Z(\xi)|^2\Big)\\
	        &\geq \varepsilon |\partial_\alpha Z|_{\inf}^2 |\alpha-\gamma|^2.
	    \end{aligned}
	\end{equation*}
	 We now see that we can take $\varepsilon>0$ explicit by using the arc-chord quantity. For all $\gamma\in\mathbb{R}^2\setminus\{0\}$, and $\xi\in\mathbb{R}^2$,  
	\begin{equation*}
	    |\partial_\alpha Z|_{\inf}\leq \liminf_{\gamma\to 0}\frac{|Z(\xi+\gamma)-Z(\xi)|}{|\gamma|}\leq|\partial_\alpha Z(\xi)\frac{\gamma}{|\gamma|}|.
	\end{equation*}
	Therefore,  
	\begin{equation*}
	   |\gamma_1\partial_{\alpha_1}Z(\xi)+\gamma_2\partial_{\alpha_2}Z(\xi)|^2\geq  |\gamma|^2|\partial_\alpha Z|_{\inf}^2. 
	\end{equation*}
	Taking $\gamma_j=\pm|\partial_{\alpha_j}Z(\xi)|^{-1}$, we conclude that
	\begin{equation*}
	    |\cos{\theta(\xi)}|\leq 1-\frac{|\pa Z|_{\inf}^2}{|\partial_{\alpha} Z(\xi)|^2}\leq 1-\frac{|\pa Z|_{\inf}^2}{\|\partial_{\alpha} Z\|_{L^\infty}^2}.
	\end{equation*}
Hence we can take 
\begin{equation}\label{varepsilon_def}
\varepsilon=(\|F(Z)\|_{L^\infty}\|\pa Z\|_{L^\infty})^{-2},    
\end{equation}
and thus we have the bound
\begin{equation}\label{bound_den}
    \begin{aligned}
        |\partial_\alpha Z(\xi)(\alpha-\gamma)|&\geq  \frac{|\alpha-\gamma|}{\|F(Z)\|_{L^\infty}^2\|\pa Z\|_{L^\infty}}.
    \end{aligned}
\end{equation}
This yields the following bound for the kernel $k(\pa Z(\xi)(\alpha-\gamma))$, 
\begin{equation}\label{kernel_bound}
    \begin{aligned}
    |k(\pa Z (\xi)(\alpha-\gamma))|\leq \frac{\|F(Z)\|_{L^\infty}^4\|\pa Z\|_{L^\infty}^2}{|\alpha-\gamma|^2}. 
    \end{aligned}
\end{equation}
We notice that
	\begin{equation}\label{den_dif}
	    \begin{aligned}
	        \frac{1}{|\pa Z(\xi)(\alpha-\gamma)|^5}-\frac{1}{|\pa Z(\xi)(\beta-\gamma)|^5}      =&\frac{(\pa Z(\xi)(\beta-\alpha))\cdot\big(\pa Z(\xi)(\alpha-\gamma+\beta-\gamma)\big) }{|\pa Z(\xi)(\alpha-\gamma)|^5|\pa Z(\xi)(\beta-\gamma)|^5}\\
	        &\times \frac{p_4(|\pa Z(\xi)(\alpha-\gamma)|^2,|\pa Z(\xi)(\beta-\gamma)|^2)}{|\pa Z(\xi)(\alpha-\gamma)|^5+|\pa Z(\xi)(\beta-\gamma)|^5},
	    \end{aligned}
	\end{equation}
with $p_m$ a homogeneous polynomial of degree $m$. Thanks to \eqref{distances} and the cut-off function, in $I_{2,2}$ it holds that $\frac12|\alpha-\gamma|\leq |\beta-\gamma|\leq\frac32|\alpha-\gamma|$. 
Together with \eqref{bound_den} and \eqref{distances}, this gives that
\begin{equation}\label{kernel_dif}
    \begin{aligned}
    |k(\pa Z(\xi)(\alpha\!-\!\gamma))\!-\!k(\pa Z(\xi)(\beta\!-\!\gamma))|\!\leq\! C\frac{\|F(Z)\|_{L^\infty}^7\|\pa Z\|_{L^\infty}^4}{|\alpha-\gamma|^3}|h|.
    \end{aligned}
\end{equation}
Then, the term $I_{2,2}$ can be estimated directly since the kernel is not singular in its domain,
	\begin{equation*}
		\begin{aligned}
			|I_{2,2}|\leq C\|g\|_{L^\infty}\|\pa Z\|_{L^\infty}^5\|F(Z)\|_{L^\infty}^8\frac{|h|}{\eta}.
		\end{aligned}
	\end{equation*} 
The support of $\chi$ allows to rewrite $I_{2,1}$ as follows
$$
I_{2,1}=\int_{\R^2}\Big(k(\pa Z(\xi)(\alpha-\gamma))-k(\pa Z(\xi)(\beta-\gamma))\Big)\chi(\gamma)g(\gamma)\minspace d\gamma.
$$
The classical splitting to show that the singular integral goes from $\dot{C}^{\sigma}$ to $\dot{C}^{\sigma}$ (see Lemma 4.6 in \cite{MajdaBertozzi2002} for example) provides the desired bound:
$$
I_{2,1}\leq C P(\|F(Z)\|_{L^\infty}\!+\!\|\pa Z\|_{L^\infty}) \|g\|_{C^\sigma}|h|^\sigma.
$$
Combining the  bounds above, we have obtained that
\begin{equation}\label{I2bound_boundary}
\begin{aligned}
|I_2|\leq C P(\|F(Z)\|_{L^\infty}\!+\!\|\pa Z\|_{L^\infty}) \|g\|_{C^\sigma}|h|^\sigma.
\end{aligned}
\end{equation}
We proceed to estimate $I_1$ \eqref{Isplit_boundary}.
	We split it as follows:
	\begin{equation}\label{I1split_boundary}
		I_1=I_{1,1}+I_{1,2}+I_{1,3}+I_{1,4},
	\end{equation}
	where
	\begin{equation*}
		\begin{aligned}
			I_{1,1}&=\!\int_{Z^{-1}(A_\eta)}\!\!d\gamma\minspace g(\gamma)\bigg(\frac{(Z_1(\alpha)\!-\!Z_1(\gamma)\!-\!\pa Z_1(\xi)\cdot(\alpha\!-\!\gamma))(Z_2(\alpha)\!-\!Z_2(\gamma))(Z_3(\alpha)\!-\!Z_3(\gamma))}{|Z(\alpha)\!-\!Z(\gamma)|^5}\\
			&\hspace{2cm}-\frac{(Z_1(\beta)-Z_1(\gamma)-\pa Z_1(\xi)\cdot(\beta-\gamma))(Z_2(\beta)-Z_2(\gamma))(Z_3(\beta)-Z_3(\gamma))}{|Z(\beta)-Z(\gamma)|^5}\bigg),
		\end{aligned}
	\end{equation*}
	\begin{equation*}
		\begin{aligned}
			I_{1,2}&=\!\int_{Z^{-1}(A_\eta)}\!\!d\gamma\minspace g(\gamma)\bigg(\frac{\pa Z_1(\xi)\cdot(\alpha-\gamma)(Z_2(\alpha)\!-\!Z_2(\gamma)\!-\!\pa Z_2(\xi)\cdot(\alpha\!-\!\gamma))(Z_3(\alpha)\!-\!Z_3(\gamma))}{|Z(\alpha)\!-\!Z(\gamma)|^5}\\
			&\hspace{1.8cm}-\frac{\pa Z_1(\xi)\cdot(\beta-\gamma)(Z_2(\beta)-Z_2(\gamma)-\pa Z_2(\xi)\cdot(\beta-\gamma))(Z_3(\beta)-Z_3(\gamma))}{|Z(\beta)-Z(\gamma)|^5}\bigg),
		\end{aligned}
	\end{equation*}
	\begin{equation*}
		\begin{aligned}			I_{1,3}&=\!\int_{Z^{-1}(A_\eta)}\!\!d\gamma\minspace g(\gamma)\bigg(\prod_{i=1}^2\pa Z_i(\xi)\cdot(\alpha-\gamma)\frac{(Z_3(\alpha)\!-\!Z_3(\gamma)-\pa Z_3(\xi)\cdot(\alpha-\gamma))}{|Z(\alpha)\!-\!Z(\gamma)|^5}\\
			&\hspace{3.4cm}-\prod_{i=1}^2\pa Z_i(\xi)\cdot(\beta-\gamma)\frac{(Z_3(\beta)-Z_3(\gamma)-\pa Z_3(\xi)\cdot(\beta-\gamma))}{|Z(\beta)-Z(\gamma)|^5}\bigg),
		\end{aligned}
	\end{equation*}
	\begin{equation*}
		\begin{aligned}
			I_{1,4}&=\!\int_{Z^{-1}(A_\eta)}\!\!d\gamma\minspace g(\gamma)\Bigg(\prod_{i=1}^3\pa Z_i(\xi)\cdot(\alpha-\gamma)
			\Big(\frac{1}{|Z(\alpha)\!-\!Z(\gamma)|^5}-\frac{1}{|\pa Z(\xi)(\alpha-\gamma)|^5}\Big)
			\\
			&\hspace{3.4cm}-\prod_{i=1}^3\pa Z_i(\xi)\cdot(\beta-\gamma)\Big(\frac{1}{|Z(\beta)-Z(\gamma)|^5}-\frac{1}{|\pa Z(\xi)(\beta-\gamma)|^5}\Big)\Bigg).
		\end{aligned}
	\end{equation*}
	Using the following sets
	\begin{equation}\label{U1U2}
	U_1=Z^{-1}(A_\eta)\cap \{|\alpha-\gamma|<2|\alpha-\beta|\},\quad U_2=Z^{-1}(A_\eta)\cap \{|\alpha-\gamma|\geq2|\alpha-\beta|\},
	\end{equation}
we decompose $I_{1,1}$ further,
	\begin{equation}\label{I11split_boundary}
		I_{1,1}=J_1+J_2+J_3+J_4+J_5,
	\end{equation}
	with
	\begin{equation*}
		\begin{aligned}
			J_1&=\int_{U_1}d\gamma\minspace g(\gamma)\bigg( \frac{(Z_1(\alpha)\!-\!Z_1(\gamma)\!-\!\pa Z_1(\xi)\cdot(\alpha\!-\!\gamma))(Z_2(\alpha)\!-\!Z_2(\gamma))(Z_3(\alpha)\!-\!Z_3(\gamma))}{|Z(\alpha)\!-\!Z(\gamma)|^5}\\
			&\hspace{2.2cm}-\frac{(Z_1(\beta)-Z_1(\gamma)-\pa Z_1(\xi)\cdot(\beta-\gamma))(Z_2(\beta)-Z_2(\gamma))(Z_3(\beta)-Z_3(\gamma))}{|Z(\beta)-Z(\gamma)|^5}\bigg)\\
			&=J_{1,1}+J_{1,2},
		\end{aligned}
	\end{equation*}
	\begin{equation*}
		\begin{aligned}
			J_2&=\int_{U_2}d\gamma\minspace g(\gamma)\frac{(Z_1(\alpha)\!-\!Z_1(\beta)\!-\!\pa Z_1(\xi)\cdot(\alpha\!-\!\beta))(Z_2(\alpha)\!-\!Z_2(\gamma))(Z_3(\alpha)\!-\!Z_3(\gamma))}{|Z(\alpha)\!-\!Z(\gamma)|^5},
		\end{aligned}
	\end{equation*}
	\begin{equation*}
		\begin{aligned}
			J_3&=\int_{U_2}d\gamma\minspace g(\gamma)\frac{(Z_1(\beta)\!-\!Z_1(\gamma)\!-\!\pa Z_1(\xi)\cdot(\beta\!-\!\gamma))(Z_2(\alpha)\!-\!Z_2(\beta))(Z_3(\alpha)\!-\!Z_3(\gamma))}{|Z(\alpha)\!-\!Z(\gamma)|^5},
		\end{aligned}
	\end{equation*}
	\begin{equation*}
		\begin{aligned}
			J_4&=\int_{U_2}d\gamma\minspace g(\gamma)\frac{(Z_1(\beta)\!-\!Z_1(\gamma)\!-\!\pa Z_1(\xi)\cdot(\beta\!-\!\gamma))(Z_2(\beta)\!-\!Z_2(\gamma))(Z_3(\alpha)\!-\!Z_3(\beta))}{|Z(\alpha)\!-\!Z(\gamma)|^5},
		\end{aligned}
	\end{equation*}
	\begin{equation*}
		\begin{aligned}
			J_5&=\int_{U_2}d\gamma\minspace g(\gamma)(Z_1(\beta)\!-\!Z_1(\gamma)\!-\!\pa Z_1(\xi)\cdot(\beta\!-\!\gamma))(Z_2(\beta)\!-\!Z_2(\gamma))(Z_3(\beta)\!-\!Z_3(\gamma))\\
			&\hspace{3cm}\times \Big(
			\frac{1}{|Z(\alpha)\!-\!Z(\gamma)|^5}-\frac{1}{|Z(\beta)\!-\!Z(\gamma)|^5}\Big).
		\end{aligned}
	\end{equation*}
	Then, we have that
	\begin{equation*}
	    \begin{aligned}
	        J_{1,1}&=\int_{U_1}d\gamma\minspace g(\gamma) \frac{(Z_1(\alpha)\!-\!Z_1(\gamma)\!-\!\pa Z_1(\alpha)\cdot(\alpha\!-\!\gamma))(Z_2(\alpha)\!-\!Z_2(\gamma))(Z_3(\alpha)\!-\!Z_3(\gamma))}{|Z(\alpha)\!-\!Z(\gamma)|^5}\\
	        &\quad+\int_{U_1}d\gamma\minspace g(\gamma) \frac{(\pa Z_1(\alpha)\!-\!\pa Z_1(\xi))\cdot(\alpha\!-\!\gamma)(Z_2(\alpha)\!-\!Z_2(\gamma))(Z_3(\alpha)\!-\!Z_3(\gamma))}{|Z(\alpha)\!-\!Z(\gamma)|^5}\\
	        &=K_1+K_2.
	    \end{aligned}
	\end{equation*}
The estimate for $K_1$ follows immediately
\begin{equation*}
    \begin{aligned}
        |K_1|&\leq C\|g\|_{L^\infty}\|\pa Z\|_{\dot{C}^\sigma}\|F(Z)\|_{L^\infty}^{6+\sigma}\|\pa Z\|_{L^\infty}^3|h|^\sigma.
    \end{aligned}
\end{equation*}
Since $U_1$ is a ball, we can write $K_2$ as follows
	\begin{equation*}
		\begin{aligned}
			K_2&=(\pa Z_1(\alpha)-\pa Z_1(\xi))\cdot\int_{U_1}d\gamma\minspace (\alpha-\gamma)\bigg(g(\gamma)\frac{(Z_2(\alpha)-Z_2(\gamma))(Z_3(\alpha)-Z_3(\gamma))}{|Z(\alpha)-Z(\gamma)|^5}\\
			&\hspace{6.5cm}-g(\alpha)\frac{\pa Z_2(\alpha)\cdot(\alpha-\gamma)\pa Z_3(\alpha)\cdot(\alpha-\gamma)}{|\pa Z(\alpha)(\alpha-\gamma)|^5}\bigg).
		\end{aligned}
	\end{equation*}
Similarly as in \eqref{den_dif}, the difference between the denominators is given by
\begin{equation}\label{den_dif2}
    \begin{aligned}
    \frac{(Z_2(\alpha)-Z_2(\gamma))(Z_3(\alpha)-Z_3(\gamma))}{|Z(\alpha)-Z(\gamma)|^5}&-\frac{\pa Z_2(\alpha)\cdot(\alpha-\gamma)\pa Z_3(\alpha)\cdot(\alpha-\gamma)}{|\pa Z(\alpha)(\alpha-\gamma)|^5}\\
    &\hspace{-4cm}=\frac{(Z_2(\alpha)-Z_2(\gamma)-\pa Z_2(\alpha)\cdot(\alpha-\gamma))\cdot(Z_2(\alpha)-Z_2(\gamma)+\pa Z_2(\alpha)\cdot(\alpha-\gamma))}{|Z(\alpha)-Z(\gamma)|^5|\pa Z(\alpha)(\alpha-\gamma)|^5}\\
    &\hspace{-4cm}\quad\times\frac{p_4(|Z(\alpha)-Z(\gamma)|^2,|\pa Z(\alpha)(\alpha-\gamma)|^2)}{|Z(\alpha)-Z(\gamma)|^5+|\pa Z(\alpha)(\alpha-\gamma)|^5},
    \end{aligned}
\end{equation}
thus
\begin{equation*}
    \begin{aligned}
        |K_2|&\leq C\big(\|g\|_{\dot{C}^\sigma}\|F(Z)\|_{L^\infty}^{6+\sigma}\|\pa Z\|_{L^\infty}^4+ \|g\|_{L^\infty}\|\pa Z\|_{\dot{C}^{\sigma}}\|F(Z)\|_{L^\infty}^{8+\sigma}\|\pa Z\|_{L^\infty}^5\big)|h|^\sigma.
    \end{aligned}
\end{equation*}
It is clear that $J_{1,2}$ satisfies the same bound, hence
\begin{equation*}
    \begin{aligned}
        |J_1|&\leq  C P(\|F(Z)\|_{L^\infty}\!+\!\|\pa Z\|_{L^\infty})\big(\|g\|_{\dot{C}^\sigma}\!+\! \|g\|_{L^\infty}\|\pa Z\|_{\dot{C}^{\sigma}}\big)|h|^\sigma.
    \end{aligned}
\end{equation*}
By writing
\begin{equation*}
\begin{aligned}
    Z_1(\alpha)-Z_1(\beta)-\pa Z_1(\xi)\cdot(\alpha-\beta)&=    Z_1(\alpha)-Z_1(\beta)-\pa Z_1(\alpha)\cdot(\alpha-\beta)\\
    &\quad+(\pa Z_1(\alpha)-\pa Z_1(\xi))\cdot(\alpha-\beta),
\end{aligned}
\end{equation*}
the term $J_2$ is bounded directly by
\begin{equation*}
    \begin{aligned}
        |J_2|&\leq C\|g\|_{L^\infty}\|\pa Z\|_{\dot{C}^\sigma}\|F(Z)\|_{L^\infty}^{6+\sigma}\|\pa Z\|_{L^\infty}^3|h|^\sigma.
    \end{aligned}
\end{equation*}
 Taking into account that on $U_2$ it holds that $\frac12|\alpha-\gamma|\leq |\beta-\gamma|\leq\frac32|\alpha-\gamma|$, the bounds for $J_3$, $J_4$, and $J_5$ follows
 \begin{equation*}
     \begin{aligned}
         |J_3|+|J_4|\leq C\|g\|_{L^\infty}\|\pa Z\|_{\dot{C}^{\sigma}}\|F(Z)\|_{L^\infty}^{6+\sigma}\|\pa Z\|_{L^\infty}^3|h|^\sigma,\\
        |J_5| \leq C\|g\|_{L^\infty}\|\pa Z\|_{\dot{C}^{\sigma}}\|F(Z)\|_{L^\infty}^{8+\sigma}\|\pa Z\|_{L^\infty}^5|h|^\sigma.
     \end{aligned}
 \end{equation*}
Inserting back in \eqref{I11split_boundary} the bounds for $J_3$-$J_7$, we obtain
\begin{equation}\label{I11boundarybound}
    \begin{aligned}
        |I_{1,1}|\leq  C P(\|F(Z)\|_{L^\infty}\!+\!\|\pa Z\|_{L^\infty})\big(\|g\|_{\dot{C}^\sigma}\!+\! \|g\|_{L^\infty}\|\pa Z\|_{\dot{C}^{\sigma}}\big)|h|^\sigma.
    \end{aligned}
\end{equation}
Now, recalling the splitting for $I_1$ \eqref{I1split_boundary}, it is clear that the bound above works as well for $I_{1,2}$, $I_{1,3}$, and $I_{1,4}$, hence it is valid for $I_1$. Combining it with the bound for $I_2$ \eqref{I2bound_boundary} in \eqref{Isplit_boundary} and recalling \eqref{holdersplit_boundary} and \eqref{gbound}, we finally have that
\begin{equation}\label{bound_boundary}
\begin{aligned}
	|S(f)(x)\!-\!S(f)(x+h)|\leq C(1\!+\!|\partial D|)P(\|F(Z)\|_{L^\infty}\!+\!\|\pa Z\|_{L^\infty})\Big(\|f\|_{C^{\sigma}}\!+\!\|f\|_{L^\infty}\|\pa Z\|_{\dot{C}^{\sigma}}\Big)|h|^\sigma.
\end{aligned}
\end{equation}
	
\vspace{0.2cm}

	\subsection{Regularity near the boundary.}\label{near_boundary}
	
	Consider two points $x\in \overline{D}$ and $x+h\in  D$ (or analogous situation in $\overline{\R^3\smallsetminus  D}$, $\R^3\smallsetminus  D$) and, without loss of generality, suppose that 
	\begin{equation}\label{distxh}
		d(x+h,\partial D)\geq d(x,\partial D)=\delta\geq0.  
	\end{equation}
	We can write
	\begin{equation}\label{xxh}
		\begin{aligned}
			x&=Z(\alpha)+\delta \tilde{N}(\alpha),\\ x+h&=Z(\alpha)+(h_n+\delta)\tilde{N}(\alpha)+h_{\tau_1}\partial_{\alpha_1} Z(\alpha)+h_{\tau_2}\partial_{\alpha_2} Z(\alpha).
		\end{aligned}
	\end{equation}
	We denote
	\begin{equation*}
		\tilde{N}(\alpha)=\frac{N(\alpha)}{\sqrt{|N(\alpha)|}}, \quad N(\alpha)=\partial_{\alpha_1}Z(\alpha)
		\wedge \partial_{\alpha_2}Z(\alpha),
	\end{equation*}
	and we notice that, as in \eqref{aux}, we have
	\begin{equation*}
	    \begin{aligned}
	        |N(\alpha)|\geq (\|F(Z)\|_{L^\infty}\|\pa Z\|_{L^\infty})^{-1}|\partial_{\alpha_1}Z(\alpha)||\partial_{\alpha_2}Z(\alpha)|,
	    \end{aligned}
	\end{equation*}
	hence
	\begin{equation}\label{tildeN}
	    \begin{aligned}
	        \|\pa Z\|_{L^\infty}\geq|\tilde{N}(\alpha)|\geq \|F(Z)\|_{L^\infty}^{-\frac32}\|\pa Z\|_{L^\infty}^{-\frac12}.
	    \end{aligned}
	\end{equation}
	We define the cutoffs
	\begin{equation}\label{delta_cutoff}
		\begin{aligned}
			|\delta|&\leq \frac16\Big(\frac{|\pa Z|_{\inf}}{18\|\pa Z\|_{\dot{C}^\sigma}}\Big)^{\frac{1}{\sigma}}\Big(\frac{|\pa Z|_{\inf}}{\|\pa Z\|_{L^\infty}}\Big)^{\frac12},\\
			|h_n|&\leq  \frac1{24}\Big(\frac{|\pa Z|_{\inf}}{18\|\pa Z\|_{\dot{C}^\sigma}}\Big)^{\frac{1}{\sigma}}\Big(\frac{|\pa Z|_{\inf}}{\|\pa Z\|_{L^\infty}}\Big)^{\frac12},\qquad
			|h_\tau|&\leq  \frac1{24}\Big(\frac{|\pa Z|_{\inf}}{18\|\pa Z\|_{\dot{C}^\sigma}}\Big)^{\frac{1}{\sigma}}\Big(\frac{|\pa Z|_{\inf}}{\|\pa Z\|_{L^\infty}}\Big)^{\frac12},
		\end{aligned}
	\end{equation}
	where $|h_\tau|^2=h_{\tau_1}^2+h_{\tau_2}^2$.
	It will be sometimes convenient to write the point $x+h$ in the following form 
	\begin{equation}\label{xh_normal}
		x+h=Z(\alpha+\lambda)+\mu\tilde{N}(\alpha+\lambda),
	\end{equation}
	where by assumption \eqref{distxh}
	\begin{equation}\label{mudelta}
		\mu|\tilde{N}(\alpha+\lambda)|\geq \delta |\tilde{N}(\alpha)|.
	\end{equation}
	We must first make sure that such a $\lambda$ and $\mu$ always exist for given $\delta$ and $h$ satisfying \eqref{delta_cutoff}.
	More specifically, we want to find $\lambda=(\lambda_1,\lambda_2)$ and $\mu$ satisfying \eqref{mudelta} solutions of \eqref{xh_normal},
	\begin{equation*}
		Z(\alpha)+(h_n+\delta)\tilde{N}(\alpha)+h_{\tau_1}\partial_{\alpha_1} Z(\alpha)+h_{\tau_2}\partial_{\alpha_2} Z(\alpha)=Z(\alpha+\lambda)+\mu\tilde{N}(\alpha+\lambda).
	\end{equation*}
	Let us denote
	\begin{equation*}
		\bm{\lambda}=(\lambda_1,\lambda_2,\mu)^T,\qquad \bm{h}=(h_{\tau_1},h_{\tau_2},h_n+\delta)^T.
	\end{equation*}
	Then, upon projecting the equations onto $\partial_{\alpha_1}Z(\alpha)$, $\partial_{\alpha_2}Z(\alpha)$, and $\tilde{N}(\alpha)$, the system reads as follows
	\begin{equation*}
		\begin{aligned}
			M(\lambda)\bm{\lambda}=\tilde{\bm{h}},
		\end{aligned}
	\end{equation*}
	where
	\begin{equation*}
		\begin{aligned}
			M(\lambda)=\begin{bmatrix}\frac{Z(\alpha+\lambda_1 e_1)-Z(\alpha)}{\lambda_1}\cdot\frac{\partial_{\alpha_1}Z(\alpha)}{|\partial_{\alpha_1}Z(\alpha)|^2} & \frac{Z(\alpha+\lambda)-Z(\alpha+\lambda_1 e_1)}{\lambda_2}\cdot\frac{\partial_{\alpha_1}Z(\alpha)}{|\partial_{\alpha_1}Z(\alpha)|^2} & \frac{\tilde{N}(\alpha+\lambda)\cdot\partial_{\alpha_1}Z(\alpha)}{|\partial_{\alpha_1}Z(\alpha)|^2}\\[2ex]
				\frac{Z(\alpha+\lambda)-Z(\alpha+\lambda_2 e_2)}{\lambda_1}\cdot\frac{\partial_{\alpha_2}Z(\alpha)}{|\partial_{\alpha_2}Z(\alpha)|^2} & \frac{Z(\alpha+\lambda_2 e_2)-Z(\alpha)}{\lambda_2}\cdot\frac{\partial_{\alpha_2}Z(\alpha)}{|\partial_{\alpha_2}Z(\alpha)|^2} & \frac{\tilde{N}(\alpha+\lambda)\cdot\partial_{\alpha_2}Z(\alpha)}{|\partial_{\alpha_2}Z(\alpha)|^2}\\[2ex]
				\frac{Z(\alpha+\lambda)-Z(\alpha+\lambda_2e_2)}{\lambda_1}\cdot\frac{\tilde{N}(\alpha)}{|\tilde{N}(\alpha)|^2} & \frac{Z(\alpha+\lambda_2e_2)-Z(\alpha)}{\lambda_2}\cdot\frac{\tilde{N}(\alpha)}{|\tilde{N}(\alpha)|^2} & \frac{\tilde{N}(\alpha+\lambda)}{|\tilde{N}(\alpha)|}\cdot\frac{\tilde{N}(\alpha)}{|\tilde{N}(\alpha)|}
			\end{bmatrix},
		\end{aligned}
	\end{equation*}
	and 
	\begin{equation*}
    \tilde{h}_1=h_{\tau_1}+h_{\tau_2}\frac{\partial_{\alpha_1}Z(\alpha)\cdot\partial_{\alpha_2}Z(\alpha)}{|\partial_{\alpha_1}Z(\alpha)|^2},\quad \tilde{h}_2=h_{\tau_2}+h_{\tau_1}\frac{\partial_{\alpha_1}Z(\alpha)\cdot\partial_{\alpha_2}Z(\alpha)}{|\partial_{\alpha_2}Z(\alpha)|^2},\quad \tilde{h}_3=h_{\tau_3}.
	\end{equation*}
	Denote the limit matrix by  $\tilde{M}(\alpha)$,
	\begin{equation*}
	    \tilde{M}=\begin{bmatrix}
	    1 & \frac{\partial_{\alpha_1}Z(\alpha)\cdot\partial_{\alpha_2}Z(\alpha)}{|\partial_{\alpha_1}Z(\alpha)|^2} & 0\\
	    \frac{\partial_{\alpha_1}Z(\alpha)\cdot\partial_{\alpha_2}Z(\alpha)}{|\partial_{\alpha_2}Z(\alpha)|^2} & 1 & 0\\
	    0 & 0 & 1
	    \end{bmatrix}.
	\end{equation*}
	As a fixed point equation, the equation reads \begin{equation}\label{fixedpointeq1}
		\bm{\lambda}=T(\bm{\lambda}):=\tilde{M}^{-1}\big((\tilde{M}-M(\lambda))\bm{\lambda}+\tilde{\bm{h}}\big).
	\end{equation}
	Since $\tilde{M}^{-1}\tilde{\bm{h}}=\bm{h}$, and taking into account \eqref{aux}, we see that 
	\begin{equation*}
		\begin{aligned}
		 |T(\bm{\lambda})|&\leq |\tilde{M}^{-1}(\tilde{M}-M(\lambda))\bm{\lambda}|+|\bm{h}|\\
		&\leq C\|\pa Z\|_{L^\infty}^2\|F(Z)\|_{L^\infty}^3\|\pa Z\|_{\dot{C}^\sigma}(1+\|\pa Z\|_{L^\infty}^{\frac12}\|F(Z)\|_{L^\infty}^{\frac12})|\bm{\lambda}|^{1+\sigma}+|\bm{h}|,
		\end{aligned}
	\end{equation*}
	thus there exists a large enough constant $C_1>0$ so that for 
	$$|\bm{h}|<\Big(C_1\|\pa Z\|_{L^\infty}^2\|F(Z)\|_{L^\infty}^3\|\pa Z\|_{\dot{C}^\sigma}(1+\|\pa Z\|_{L^\infty}^{\frac12}\|F(Z)\|_{L^\infty}^{\frac12})\Big)
	^{-\frac{1}{\sigma}},$$
	Brouwer's Fixed Point Theorem \cite{Evans} yields the existence of a solution to \eqref{fixedpointeq1} in the ball of radius $2\Big(C_1\|\pa Z\|_{L^\infty}^2\|F(Z)\|_{L^\infty}^3\|\pa Z\|_{\dot{C}^\sigma}(1+\|\pa Z\|_{L^\infty}^{\frac12}\|F(Z)\|_{L^\infty}^{\frac12})\Big)
	^{-\frac{1}{\sigma}}$ and centered at the origin, and moreover
	\begin{equation}\label{lambdabound}
		|\bm{\lambda}|\leq 2|\bm{h}|.
	\end{equation}
	Finally, the third equation in \eqref{fixedpointeq1} shows that the condition \eqref{mudelta}, i.e., \eqref{distxh}, implies that
	\begin{equation*}
		\begin{aligned}
			h_n&\geq -\delta\big|\frac{|\tilde{N}(\alpha)|-|\tilde{N}(\alpha+\lambda)|}{|\tilde{N}(\alpha+\lambda)|}\big|-3C_1^{-1}|\bm{\lambda}|,
		\end{aligned}
	\end{equation*}
	which in particular gives that, for suitable $C_1$ big enough,
	\begin{equation}\label{hn_lower}
		h_n\geq -\frac{\delta}{4}-\frac{|h_\tau|}{4}.
	\end{equation}
	\vspace{0.2cm}
	
	Next, we distinguish two cases: $|h_\tau|\leq \frac14\frac{|\pa Z|_{\inf}}{\|\pa Z\|_{L^\infty}}\delta$ and $|h_\tau|\geq \frac14\frac{|\pa Z|_{\inf}}{\|\pa Z\|_{L^\infty}}\delta$.
	
	\vspace{0.2cm}
	
	\subsubsection{{\em{\textbf{Regularity in {\em{nearly}} normal direction:}}}} Assume that \begin{equation}\label{htau_menor_delta}
		|h_\tau|\leq \frac14\frac{|\pa Z|_{\inf}}{\|\pa Z\|_{L^\infty}}\delta.
	\end{equation}
	We can write
	\begin{equation}\label{normaltangentSplit}
		\begin{aligned}
			S(f)(x+h)-S(f)(x)&=S(f)(x+h)-S(f)(Z(\alpha)+(\delta+h_n)\tilde{N}(\alpha))\\
			&\quad+S(f)(Z(\alpha)+(\delta+h_n)\tilde{N}(\alpha))-S(f)(Z(\alpha)+\delta \tilde{N}(\alpha)),
		\end{aligned}
	\end{equation}
	so that the first two terms correspond to a difference in the tangential direction and the last two to a difference in the normal direction.
	We estimate each of these terms separately.
	Note that the above splitting is valid since \eqref{hn_lower} and the assumption $|h_\tau|\leq \frac14\frac{|\pa Z|_{\inf}}{\|\pa Z\|_{L^\infty}}\delta$ guarantees that $h_n+\delta\geq \delta/2$, hence the point $Z(\alpha)+(\delta+h_n)\tilde{N}(\alpha)$ belongs to $ D$. We can thus assume in the following subsection, Case 1, that $h_n\geq0$; otherwise, interchange the roles of $\delta$ and $\delta+h_n$.
	\vspace{0.2cm}
	
	\noindent\textbf{Case 1: Normal direction.} 
	Here we consider $h_\tau=0$, i.e., we are dealing with the second difference above. We write it as follows
	\begin{equation*}
		\begin{aligned}
			S(f)(Z(\alpha)+(\delta+&h_n)\tilde{N}(\alpha))-S(f)(Z(\alpha)+\delta \tilde{N}(\alpha))\\
			&=\pv\int_{Z(B_\eta)}\Big(k(Z(\alpha)+(\delta+h_n)\tilde{N}(\alpha)-y)-k(x-y)\Big)f(y)dS(y)\\
			&\quad+\int_{\partial D \smallsetminus Z(B_\eta)}\Big(k(Z(\alpha)+(\delta+h_n)\tilde{N}(\alpha)-y)-k(x-y)\Big)f(y)dS(y)\\
			&=I+II,
		\end{aligned}
	\end{equation*}
	where 
	\begin{equation*}
		B_\eta=\{\gamma\in\mathbb{R}^2:|\alpha-\gamma|<\eta\}.
	\end{equation*}
	The second term is again more regular, 
	\begin{equation*}
		\begin{aligned}
			|II|\leq \frac{C}{\eta^3|\pa Z|_{\inf}^3}|\partial D|\|f\|_{L^\infty}|h|. 
		\end{aligned}
	\end{equation*}
	The first term is given by
	\begin{equation*}
		\begin{aligned}
			I&=\int_{B_\eta} \frac{\prod_{j=1}^3(Z_j(\alpha)+(h_n+\delta)\tilde{N}_j(\alpha)-Z_j(\gamma))}{|Z(\alpha)+(h_n+\delta)\tilde{N}(\alpha)-Z(\gamma)|^5}g(\gamma)d\gamma \\
			&\quad-\int_{B_\eta} \frac{\prod_{j=1}^3(Z_j(\alpha)+\delta\tilde{N}_j(\alpha)-Z_j(\gamma))}{|Z(\alpha)+\delta\tilde{N}(\alpha)-Z(\gamma)|^5}g(\gamma)d\gamma,
		\end{aligned}    
	\end{equation*}
	and we decompose it as follows
	\begin{equation}\label{Isplit}
		I=I_1+I_2,
	\end{equation}
	\begin{equation*}
		\begin{aligned}
			I_1&=\int_{B_\eta}\prod_{j=1}^3(Z_j(\alpha)-Z_j(\gamma))\Big(\frac{1}{|Z(\alpha)+(\delta+h_n)\tilde{N}(\alpha)-Z(\gamma)|^5}-\frac{1}{|x-Z(\gamma)|^5}\Big)g(\gamma)d\gamma,
		\end{aligned}
	\end{equation*}
	\begin{equation}
		\begin{aligned}\label{I2N}
			I_2&=\int_{B_\eta}\!\!\frac{\prod_{j=1}^3(Z_j(\alpha)-Z_j(\gamma)+(h_n+\delta)\tilde{N}_j(\alpha))-\prod_{j=1}^3(Z_j(\alpha)-Z_j(\gamma))}{|Z(\alpha)-Z(\gamma)+(h_n+\delta)\tilde{N}(\alpha)|^5}g(\gamma)d\gamma \\
			&\quad-\int_{B_\eta}\frac{\prod_{j=1}^3(Z_j(\alpha)-Z_j(\gamma)+\delta\tilde{N}_j(\alpha))-\prod_{j=1}^3(Z_j(\alpha)-Z_j(\gamma))}{|Z(\alpha)-Z(\gamma)+\delta\tilde{N}(\alpha)|^5}g(\gamma)d\gamma.
		\end{aligned}
	\end{equation}
	We split $I_1$ further:
	\begin{equation}\label{I1split}
		\begin{aligned}
			I_1&=I_{1,1}+\dots+I_{1,6},
		\end{aligned}
	\end{equation}
	where
	\begin{equation*}
		\begin{aligned}
			I_{1,1}&=\int_{B_\eta}(Z_1(\alpha)-Z_1(\gamma)-\pa Z_1(\alpha)\cdot(\alpha-\gamma))\prod_{j=1}^2(Z_j(\alpha)-Z_j(\gamma))\\
			&\hspace{1cm}\times \Big(\frac{1}{|Z(\alpha)+(\delta+h_n)\tilde{N}(\alpha)-Z(\gamma)|^5}-\frac{1}{|x-Z(\gamma)|^5}\Big)g(\gamma)d\gamma,
		\end{aligned}
	\end{equation*}
	\begin{equation*}
		\begin{aligned}
			I_{1,2}&=\int_{B_\eta}\pa Z_1(\alpha)\cdot(\alpha-\gamma)(Z_2(\alpha)-Z_2(\gamma)-\pa Z_2(\alpha)\cdot(\alpha-\gamma))(Z_3(\alpha)-Z_3(\gamma))\\
			&\hspace{1cm}\times\Big(\frac{1}{|Z(\alpha)+(\delta+h_n)\tilde{N}(\alpha)-Z(\gamma)|^5}-\frac{1}{|x-Z(\gamma)|^5}\Big)g(\gamma)d\gamma,
		\end{aligned}
	\end{equation*}
	\begin{equation*}
		\begin{aligned}
			I_{1,3}&=\int_{B_\eta}\pa Z_1(\alpha)\cdot(\alpha-\gamma)\pa Z_2(\alpha)\cdot(\alpha-\gamma)(Z_3(\alpha)-Z_3(\gamma)-\pa Z_3(\alpha-\gamma))\\
			&\hspace{1cm}\times\Big(\frac{1}{|Z(\alpha)+(\delta+h_n)\tilde{N}(\alpha)-Z(\gamma)|^5}-\frac{1}{|x-Z(\gamma)|^5}\Big)g(\gamma)d\gamma,
		\end{aligned}
	\end{equation*}
	\begin{equation*}
		\begin{aligned}
			I_{1,4}&=\int_{B_\eta}\pa Z_1(\alpha)\cdot(\alpha-\gamma)\pa Z_2(\alpha)\cdot(\alpha-\gamma)\pa Z_3(\alpha)\cdot(\alpha-\gamma)\\
			&\hspace{1cm}\times (g(\gamma)-g(\alpha))\Big(\frac{1}{|Z(\alpha)+(\delta+h_n)\tilde{N}(\alpha)-Z(\gamma)|^5}-\frac{1}{|x-Z(\gamma)|^5}\Big)d\gamma,
		\end{aligned}
	\end{equation*}
	\begin{equation*}
		\begin{aligned}
			I_{1,5}&=g(\alpha)\int_{B_\eta}\pa Z_1(\alpha)\cdot(\alpha-\gamma)\pa Z_2(\alpha)\cdot(\alpha-\gamma)\pa Z_3(\alpha)\cdot(\alpha-\gamma)\\
			&\hspace{1cm}\times \Big(\frac{1}{|Z(\alpha)+(\delta+h_n)\tilde{N}(\alpha)-Z(\gamma)|^5}-\frac{1}{|\pa Z(\alpha)(\alpha-\gamma)+(h_n+\delta)\tilde{N}(\alpha)|^5}\\
			&\hspace{1cm}\qquad+\frac{1}{|\pa Z(\alpha) (\alpha-\gamma)+\delta \tilde{N}(\alpha)|^5}-\frac{1}{|x-Z(\gamma)|^5}\Big)d\gamma,
		\end{aligned}
	\end{equation*}
	\begin{equation*}
		\begin{aligned}
			I_{1,6}&=g(\alpha)\int_{B_\eta}\pa Z_1(\alpha)\cdot(\alpha-\gamma)\pa Z_2(\alpha)\cdot(\alpha-\gamma)\pa Z_3(\alpha)\cdot(\alpha-\gamma)\\
			&\hspace{1cm}\times \Big(\frac{1}{|\pa Z(\alpha)(\alpha-\gamma)+(h_n+\delta)\tilde{N}(\alpha)|)^{5}}-\frac{1}{|\pa Z(\alpha)(\alpha-\gamma)+\delta\tilde{N}(\alpha)|^5}\Big)d\gamma,
		\end{aligned}
	\end{equation*}
	To estimate these terms we will need to bound from below the denominator
	\begin{equation*}
		D=|Z(\alpha)-Z(\gamma)+(\delta+h_n)\tilde{N}(\alpha)|^2.
	\end{equation*}
	We can write
	\begin{equation*}
		\begin{aligned}
			D&=|Z(\alpha)-Z(\gamma)|^2+(h_n+\delta)^2|N(\alpha)|+2(Z(\alpha)-Z(\gamma)-\pa Z(\alpha)(\alpha-\gamma))\cdot\tilde{N}(\alpha)(h_n+\delta)\\
			&\geq \frac{|Z(\alpha)-Z(\gamma)|^2}{|\alpha-\gamma|^2} |\alpha-\gamma|^2+(h_n+\delta)^2|\tilde{N}(\alpha)|^2-2|\alpha-\gamma|^{1+\sigma}\|\pa Z\|_{\dot{C}^\sigma}|\tilde{N}(\alpha)|(h_n+\delta).
		\end{aligned}
	\end{equation*}
    The last term satisfies that
	\begin{equation*}
		\begin{aligned}
			2|\alpha-\gamma|^{1+\sigma}\|\pa Z\|_{\dot{C}^\sigma}|\tilde{N}(\alpha)|(h_n+\delta)&\leq\frac{1-\sigma}2(h_n+\delta)^2|\tilde{N}(\alpha)|^2\\
			+ \frac{1+\sigma}2 2^{\frac{2}{1+\sigma}}&(h_n+\delta)^{\frac{2\sigma}{1+\sigma}}\|\pa Z\|_{L^\infty}^{\frac{2\sigma}{1+\sigma}}\|\pa Z\|_{\dot{C}}^{\frac{2}{1+\sigma}}|\alpha-\gamma|^2.
		\end{aligned}
	\end{equation*}
	The fact that $h_n+\delta\geq\delta/2$ and the choice of the cutoff for $\delta$ \eqref{delta_cutoff}  allow us to obtain that
	\begin{equation}\label{denominatorb2}
		\begin{aligned}
			D&\geq \frac12 |\pa Z|_{\text{inf}}^2\Big( |\alpha-\gamma|^2+(h_n+\delta)^2\Big).
		\end{aligned}
	\end{equation}
	Next, we proceed to estimate each of these terms $I_{1,i}$, $i=1,\dots, 6$. For $I_{1,1}$ we have that
	\begin{equation}\label{I11split2}
		\begin{aligned}
			|I_{1,1}|\leq C|h_n|\|g\|_{L^\infty}\|\pa Z\|_{L^\infty}^3\|\pa Z\|_{\dot{C}^{\sigma}}\int_{B_\eta}  \Big(&\frac{|\alpha-\gamma|^{3+\sigma}|x-Z(\gamma)|^{-5}}{|Z(\alpha)+(\delta+h_n)\tilde{N}(\alpha)-Z(\gamma)|}\\
			&+\frac{|\alpha-\gamma|^{3+\sigma}|x-Z(\gamma)|^{-1}}{{|Z(\alpha)+(\delta+h_n)\tilde{N}(\alpha)-Z(\gamma)|^5}}\Big)d\gamma.
		\end{aligned}
	\end{equation}
	We introduce the bound for the denominator \eqref{denominatorb2} in $I_{1,1}$ to obtain that
	\begin{equation*}
		\begin{aligned}
			|I_{1,1}|\leq C|h_n|\|g\|_{L^\infty}\|\pa Z\|_{\dot{C}^{\sigma}}\frac{\|\pa Z\|_{L^\infty}^3}{|\pa Z|_{\text{inf}}^6}\int_{B_\eta}\Big(&\frac{|\alpha-\gamma|^{3+\sigma}(|\alpha-\gamma|^2+\delta^2)^{-\frac52}}{(|\alpha-\gamma|^2+(h_n+\delta)^2)^{\frac12}}\\
			&+\frac{|\alpha-\gamma|^{3+\sigma}(|\alpha-\gamma|^2+\delta^2)^{-\frac12}}{(|\alpha-\gamma|^2+(h_n+\delta)^2)^{\frac52}}\Big)d\gamma.
		\end{aligned}
	\end{equation*}
	Changing variables $w=(\alpha-\gamma)/(h_n+\delta)$,
	\begin{equation*}
		\begin{aligned}
			|I_{1,1}|\leq C\frac{|h_n|}{(h_n+\delta)^{1-\sigma}}\|g\|_{L^\infty}\|\pa Z\|_{\dot{C}^{\sigma}}\frac{\|\pa Z\|_{L^\infty}^3}{|\pa Z|_{\text{inf}}^6}\int_{|w|\leq\frac{\eta}{h_n+\delta}}&\Big(\frac{|w|^{3+\sigma}\big(|w|^2+\big(\frac{\delta}{h_n+\delta}\big)^2\big)^{-\frac52}}{(|w|^2+1)^{\frac12}}\\
			&+\frac{|w|^{3+\sigma}\big(|w|^2+\big(\frac{\delta}{h_n+\delta}\big)^2\big)^{-\frac12}}{(|w|^2+1)^{\frac52}}\Big)dw,
		\end{aligned}
	\end{equation*}
	hence
	\begin{equation*}
		\begin{aligned}
			|I_{1,1}|&\leq C|h_n|^\sigma\|g\|_{L^\infty}\|\pa Z\|_{\dot{C}^{\sigma}}\frac{\|\pa Z\|_{L^\infty}^3}{|\pa Z|_{\text{inf}}^6}\int_{\R^2}\Big(\frac{|w|^{-2+\sigma}}{(|w|^2+1)^{\frac12}}+\frac{|w|^{2+\sigma}}{(|w|^2+1)^{\frac52}}\Big)dw,
		\end{aligned}
	\end{equation*}
	to conclude the desired bound
	\begin{equation}\label{I11bound}
		\begin{aligned}
			|I_{1,1}|&\leq C\|g\|_{L^\infty}\|\pa Z\|_{\dot{C}^{\sigma}}\frac{\|\pa Z\|_{L^\infty}^3}{|\pa Z|_{\text{inf}}^6}|h_n|^\sigma.
		\end{aligned}
	\end{equation}
	The terms $I_{1,2}$ and $I_{1,3}$ are bounded analogously:
	\begin{equation}\label{I12I13bound}
		\begin{aligned}
			|I_{1,2}|+|I_{1,3}|&\leq C\|g\|_{L^\infty}\|\pa Z\|_{\dot{C}^{\sigma}}\frac{\|\pa Z\|_{L^\infty}^3}{|\pa Z|_{\text{inf}}^6}|h_n|^\sigma.
		\end{aligned}
	\end{equation}
	In $I_{1,4}$ the same approach yields
	\begin{equation}\label{I14bound2}
		|I_{1,4}|\leq C\|g\|_{\dot{C}^\sigma}\frac{\|\pa Z\|_{L^\infty}^4}{|\pa Z|_{\text{inf}}^6}|h_n|^\sigma.
	\end{equation}
	We deal with $I_{1,5}$ \eqref{I1split}. Let us denote
	\begin{equation}\label{uhvh}
		\begin{aligned}
			u_h&=|Z(\alpha)-Z(\gamma)+(h_n+\delta)\tilde{N}(\alpha)|,\\
			v_h&=|\pa Z(\alpha)(\alpha-\gamma)+(h_n+\delta)\tilde{N}(\alpha)|,
		\end{aligned}
	\end{equation}
	and 
	\begin{equation*}
		\begin{aligned}
			\frac{1}{u_h^5}-\frac{1}{u_0^5}=G(u_h,u_0)(u_0^2-u_h^2),
		\end{aligned}
	\end{equation*}
	where
	\begin{equation}\label{Guhvh}
		G(u_h,u_0)=\frac1{u_h+u_0}\big(\frac1{u_h^5u_0}+\frac1{u_h^4u_0^2}+\frac1{u_h^3u_0^3}+\frac1{u_h^2u_0^4}+\frac1{u_hu_0^5}\big).
	\end{equation}
	We notice that $u_h^2=D$ for which we have the lower bound \eqref{denominatorb2}. We also need a lower bound for $v_h$:
	\begin{equation*}
	    \begin{aligned}
	   |v_h|^2&=|\pa Z(\alpha)(\alpha-\gamma)|^2+(h_n+\delta)^2|\tilde{N}(\alpha)|^2\\
	   &\geq |\pa Z|_{\inf}^2\frac{|\pa Z|_{\inf}^2}{\|\pa Z\|_{L^\infty}^2}\big(|\alpha-\gamma|^2+(h_n+\delta)^2\big),
	    \end{aligned}
	\end{equation*}
	where we have used \eqref{aux}, \eqref{varepsilon_def}, and \eqref{tildeN}. Notice that the following is a common lower bound for $u_h$ and $v_h$,
		\begin{equation}\label{uhvh_lower}
	    \begin{aligned}
	   |u_h|^2,|v_h|^2&\geq \frac12|\pa Z|_{\inf}^2\frac{|\pa Z|_{\inf}^2}{\|\pa Z\|_{L^\infty}^2}\big(|\alpha-\gamma|^2+(h_n+\delta)^2\big).
	    \end{aligned}
	\end{equation}
	We have that
	\begin{equation}\label{uhmvh}
		\begin{aligned}
			u_0^2-u_h^2&=-h_n(h_n+2\delta)|N(\alpha)|-2(Z(\alpha)-Z(\gamma))\cdot\tilde{N}(\alpha)h_n,\\ v_0^2-v_h^2&=-h_n(h_n+2\delta)|N(\alpha)| .
		\end{aligned}
	\end{equation}
	The term $I_{1,5}$ \eqref{I1split} is then written as follows
	\begin{equation}
		\begin{aligned}\label{I15}
			I_{1,5}&=g(\alpha)\int_{B_\eta}\pa Z_1(\alpha)\cdot(\alpha-\gamma)\pa Z_2(\alpha)\cdot(\alpha-\gamma)\pa Z_3(\alpha)\cdot(\alpha-\gamma)\\
			&\hspace{2.5cm}\times \Big(\frac{1}{u_h^5}-\frac{1}{u_0^5}-\big(\frac{1}{v_h^5}-\frac{1}{v_0^5}\big)\Big)d\gamma,
		\end{aligned}
	\end{equation}
	and we split it further
	\begin{equation*}
		I_{1,5}=J_1+J_2,
	\end{equation*}
	where
	\begin{equation*}
		\begin{aligned}
			J_1&=g(\alpha)\int_{B_\eta}\pa Z_1(\alpha)\cdot(\alpha-\gamma)\pa Z_2(\alpha)\cdot(\alpha-\gamma)\pa Z_3(\alpha)\cdot(\alpha-\gamma)\\
			&\hspace{2cm}\times G(u_h,u_0)\big(u_0^2-u_h^2-(v_0^2-v_h^2)\big)d\gamma,
		\end{aligned}
	\end{equation*}
	\begin{equation*}
		\begin{aligned}
			J_2&=g(\alpha)\int_{B_\eta}\pa Z_1(\alpha)\cdot(\alpha-\gamma)\pa Z_2(\alpha)\cdot(\alpha-\gamma)\pa Z_3(\alpha)\cdot(\alpha-\gamma)\\
			&\hspace{2cm}\times \big(G(u_h,u_0)-G(v_h,v_0)\big)(v_0^2-v_h^2)d\gamma.
		\end{aligned}
	\end{equation*}
	Substituting \eqref{uhmvh},
	\begin{equation*}
		\begin{aligned}
			J_{1}=-2h_n\minspace g(\alpha)\int_{B_\eta}&\pa Z_1(\alpha)\cdot(\alpha-\gamma)\pa Z_2(\alpha)\cdot(\alpha-\gamma)\pa Z_3(\alpha)\cdot(\alpha-\gamma)\\
			&\times G(u_h,u_0)
			(Z(\alpha)-Z(\gamma))\cdot\tilde{N}(\alpha) d\gamma.
		\end{aligned}
	\end{equation*}
	The extra cancellation  
	\begin{equation*}
		(Z(\alpha)-Z(\gamma))\cdot\tilde{N}(\alpha)=(Z(\alpha)-Z(\gamma)-\pa Z(\alpha)(\alpha-\gamma))\cdot\tilde{N}(\alpha),
	\end{equation*}
	and recalling that $u_h^2=D$, we introduce  the lower bound for $D$ in \eqref{denominatorb2} to obtain that
	\begin{equation*}
		\begin{aligned}
			|J_{1}|&\leq C \minspace \|g\|_{L^\infty}\frac{\|\pa Z\|_{L^\infty}^4}{|\pa Z|_{\text{inf}}^7}\|\pa Z\|_{\dot{C}^\sigma}|h_n| \int_{|\alpha-\gamma|\leq \eta} \Big(\frac{|\alpha-\gamma|^{4+\sigma}(|\alpha-\gamma|^2+\delta^2)^{-\frac12}}{(|\alpha-\gamma|^2+(h_n+\delta)^2)^{3}} \\
			&\qquad\qquad\qquad\qquad\qquad\qquad\qquad\qquad\qquad+\frac{|\alpha-\gamma|^{4+\sigma}(|\alpha-\gamma|^2+\delta^2)^{-\frac52}}{|\alpha-\gamma|^2+(h_n+\delta)^2}\Big)d\gamma\\
			&\leq  C \minspace \|g\|_{L^\infty}\frac{\|\pa Z\|_{L^\infty}^4}{|\pa Z|_{\text{inf}}^7}\|\pa Z\|_{\dot{C}^\sigma}|h_n|^\sigma.
		\end{aligned}
	\end{equation*}
	We proceed to estimate $J_2$. We further split this term
	$$
	J_{2}=\sum_{k=1}^6J_{2,k}
	$$
	where 
	\begin{equation*}
		\begin{aligned}
			J_{2,k}&=g(\alpha)\int_{B_\eta}\pa Z_1(\alpha)\cdot(\alpha-\gamma)\pa Z_2(\alpha)\cdot(\alpha-\gamma)\pa Z_3(\alpha)\cdot(\alpha-\gamma)\\
			&\hspace{2cm}\times \big(\frac1{u_h^{6-k}u_0^k}-\frac1{v_h^{6-k}v_0^k}\big)\frac{v_0^2-v_h^2}{u_h+u_0}d\gamma,
		\end{aligned}
	\end{equation*}
	for $1\leq k\leq 5$, and
	\begin{equation*}
		\begin{aligned}
			J_{2,6}&=g(\alpha)\int_{B_\eta}\pa Z_1(\alpha)\cdot(\alpha-\gamma)\pa Z_2(\alpha)\cdot(\alpha-\gamma)\pa Z_3(\alpha)\cdot(\alpha-\gamma)(v_0^2-v_h^2)\\
			&\hspace{2cm}\times \big(\frac1{v_h^5v_0}+\frac1{v_h^4v_0^2}+\frac1{v_h^3v_0^3}+\frac1{v_h^2v_0^4}+\frac1{v_hv_0^5}\big)\big(\frac{1}{u_h+u_0}-\frac1{v_h+v_0}\big)d\gamma.
		\end{aligned}
	\end{equation*}
	To control $J_{2,1}$ a further splitting is given:
	$$
	J_{2,1}=\sum_{l=1}^6K_l
	$$
	where 
	\begin{equation*}
		\begin{aligned}
			K_{1}&=g(\alpha)\int_{B_\eta}\pa Z_1(\alpha)\cdot(\alpha-\gamma)\pa Z_2(\alpha)\cdot(\alpha-\gamma)\pa Z_3(\alpha)\cdot(\alpha-\gamma)\\
			&\hspace{2cm}\times \frac1{u_h^5} \big(\frac1{u_0}-\frac1{v_0}\big)\frac{v_0^2-v_h^2}{u_h+u_0}d\gamma,
		\end{aligned}
	\end{equation*}
	and
	\begin{equation*}
		\begin{aligned}
			K_{l}&=g(\alpha)\int_{B_\eta}\pa Z_1(\alpha)\cdot(\alpha-\gamma)\pa Z_2(\alpha)\cdot(\alpha-\gamma)\pa Z_3(\alpha)\cdot(\alpha-\gamma)\\
			&\hspace{2cm}\times \frac1{v_0u_h^{6-l}v_h^{l-2}}\big(\frac1{u_h}-\frac1{v_h}\big)\frac{v_0^2-v_h^2}{u_h+u_0}d\gamma
		\end{aligned}.
	\end{equation*}
	for $2\leq l\leq 6$. Estimate
	\begin{equation}\label{v0-u0}
		|v_h-u_h|\leq \|\partial_\alpha Z\|_{\dot{C}^{\sigma}}|\alpha-\gamma|^{1+\sigma}  
	\end{equation}
	and the lower bound \eqref{uhvh_lower} allows us to get
	\begin{equation*}
		\begin{aligned}
			|K_1|&\leq C\|g\|_{L^\infty}\frac{\|\partial_\alpha Z\|_{L^\infty}^{13}}{\|\partial_\alpha Z\|^{16}_{\text{inf}}}\|\partial_\alpha Z\|_{\dot{C}^{\sigma}}|h_n|\int_{|\alpha-\gamma|\leq \eta} \frac{|\alpha-\gamma|^{2+\sigma}((h_n+\delta)+\delta)d\gamma}{(|\alpha-\gamma|^2+(h_n+\delta)^2)^{\frac52}((h_n+\delta)+\delta)} \\
			&\leq  C \minspace \|g\|_{L^\infty}\frac{\|\pa Z\|_{L^\infty}^{13}}{|\pa Z|_{\text{inf}}^{16}}\|\pa Z\|_{\dot{C}^\sigma}|h_n|^\sigma.
		\end{aligned}
	\end{equation*}
	Using \eqref{v0-u0}, an analogous bound follows for the rest of $K_l$. It yields the desired control for the term $J_{2,1}$. 
	
	The terms  $J_{2,k}$, $k=2,...,5$, are controlled in a similar manner to $J_{2,1}$. We show some detail in the most singular one: $J_{2,5}$. It is decomposed by
	$$
	J_{2,5}=K_{r}+K_{s}
	$$
	where $K_{r}$ is a representative term, and in $K_s$ we collect similar integrals (they can be handled as before). The $K_{r}$ integral is given by
	\begin{equation*}
		\begin{aligned}
			K_{r}&=g(\alpha)\int_{B_\eta}\pa Z_1(\alpha)\cdot(\alpha-\gamma)\pa Z_2(\alpha)\cdot(\alpha-\gamma)\pa Z_3(\alpha)\cdot(\alpha-\gamma)\\
			&\hspace{2cm}\times \frac1{v_hu_0^4}(\frac1{u_0}-\frac1{v_0}\big)\frac{v_0^2-v_h^2}{u_h+u_0}d\gamma.
		\end{aligned}
	\end{equation*}
	Then, the following bound is obtained
	\begin{equation*}
		\begin{aligned}
			|K_{r}|&\leq C\|g\|_{L^\infty}\frac{\|\partial_\alpha Z\|_{L^\infty}^{13}}{\|\partial_\alpha Z\|^{16}_{\text{inf}}}\|\partial_\alpha Z\|_{\dot{C}^{\sigma}}|h_n|\int_{|\alpha-\gamma|\leq \eta} \frac{|\alpha-\gamma|^{4+\sigma}(|\alpha-\gamma|^2+\delta^2)^{-3}}{(|\alpha-\gamma|^2+(h_n+\delta)^2)^{\frac12}} d\gamma\\
			&\leq  C \minspace \|g\|_{L^\infty}\frac{\|\pa Z\|_{L^\infty}^{13}}{|\pa Z|_{\text{inf}}^{16}}\|\pa Z\|_{\dot{C}^\sigma}|h_n|^\sigma.
		\end{aligned}
	\end{equation*}
	The rest of terms in $K_s$ are estimated analogously. Then, it yields the desired control for $J_{2,5}$.
	
	Next, we consider $J_{2,6}$ which is estimated using \eqref{v0-u0}:
	\begin{equation*}
		\begin{aligned}
			|J_{2,6}|&\leq C\|g\|_{L^\infty}\frac{\|\partial_\alpha Z\|_{L^\infty}^{13}}{\|\partial_\alpha Z\|^{16}_{\text{inf}}}\|\partial_\alpha Z\|_{\dot{C}^{\sigma}}|h_n|\int_{|\alpha-\gamma|\leq \eta} \Big(\frac{|\alpha-\gamma|^{4+\sigma}(|\alpha-\gamma|^2+\delta^2)^{-\frac12}}{(|\alpha-\gamma|^2+(h_n+\delta)^2)^{3}} \\
			&\qquad\qquad\qquad\qquad\qquad\qquad\qquad\qquad\qquad+\frac{|\alpha-\gamma|^{4+\sigma}(|\alpha-\gamma|^2+\delta^2)^{-\frac52}}{|\alpha-\gamma|^2+(h_n+\delta)^2}\Big)d\gamma\\
			&\leq  C \minspace \|g\|_{L^\infty}\frac{\|\pa Z\|_{L^\infty}^{13}}{|\pa Z|_{\text{inf}}^{16}}\|\pa Z\|_{\dot{C}^\sigma}|h_n|^\sigma.
		\end{aligned}
	\end{equation*} 
	We are done with $J_2$ and therefore with $I_{1,5}$,
	\begin{equation}\label{I15bound}
	    \begin{aligned}
	        |I_{1,5}|&\leq C \minspace \|g\|_{L^\infty}\|\pa Z\|_{L^\infty}^4\|F(Z)\|_{L^\infty}^7(1+\|\pa Z\|_{L^\infty}^9\|F(Z)\|_{L^\infty}^9)\|\pa Z\|_{\dot{C}^\sigma}|h_n|^\sigma.
	    \end{aligned}
	\end{equation}
	Finally, because $|\pa Z(\alpha)(\alpha-\gamma)+(h_n+\delta)\tilde{N}(\alpha)|^2=|\pa Z(\alpha)(\alpha-\gamma)|^2+(h_n+\delta)^2|\tilde{N}(\alpha)|^2$, the integral in $I_{1,6}$ \eqref{I1split} is odd and thus vanishes. 
	Hence, recalling \eqref{I1split} and the bounds \eqref{I11bound}-\eqref{I14bound2}, \eqref{I15bound}, the $I_1$ integral is estimated,
	\begin{equation}\label{I1bound}
	    \begin{aligned}
	        |I_1|&\leq C \minspace \|g\|_{L^\infty}\|\pa Z\|_{L^\infty}^3\|F(Z)\|_{L^\infty}^6(1+\|\pa Z\|_{L^\infty}^{10}\|F(Z)\|_{L^\infty}^{10})\|\pa Z\|_{\dot{C}^\sigma}|h_n|^\sigma\\
	        &\quad+C \minspace \|g\|_{\dot{C}^\sigma}\|\pa Z\|_{L^\infty}^4\|F(Z)\|_{L^\infty}^6|h_n|^\sigma.
	    \end{aligned}
	\end{equation}
	It remains the control of $I_2$. In order to estimate the term $I_2$ a further decomposition is done in \eqref{I2N}:
	\begin{equation}\label{I2split}
		I_2=\sum_{j=1}^8I_{2,j}.
	\end{equation}
	They are given expanding the products in the numerator, gathering in one term the subtraction between one integral with $(h_n+\delta)$ and its corresponding integral with only $\delta$. Here we show how to deal with three of them, as the rest of the estimates follow similarly. They are given by
	\begin{equation*}
		\begin{aligned}
			I_{2,1}&=\int_{B_\eta}\!\!\frac{(Z_1(\alpha)-Z_1(\gamma))(Z_2(\alpha)-Z_2(\gamma))(h_n+\delta)\tilde{N}_3(\alpha)}{|Z(\alpha)-Z(\gamma)+(h_n+\delta)\tilde{N}(\alpha)|^5}g(\gamma)d\gamma \\
			&\quad-\int_{B_\eta}\frac{(Z_1(\alpha)-Z_1(\gamma))(Z_2(\alpha)-Z_2(\gamma))\delta\tilde{N}_3(\alpha)}{|Z(\alpha)-Z(\gamma)+\delta\tilde{N}(\alpha)|^5}g(\gamma)d\gamma,
		\end{aligned}
	\end{equation*}
	\begin{equation*}
		\begin{aligned}
			I_{2,2}&=\int_{B_\eta}\!\!\frac{(Z_1(\alpha)-Z_1(\gamma))(h_n+\delta)^2\tilde{N}_2(\alpha)\tilde{N}_3(\alpha)}{|Z(\alpha)-Z(\gamma)+(h_n+\delta)\tilde{N}(\alpha)|^5}g(\gamma)d\gamma \\
			&\quad-\int_{B_\eta}\frac{(Z_1(\alpha)-Z_1(\gamma))\delta^2\tilde{N}_2(\alpha)\tilde{N}_3(\alpha)}{|Z(\alpha)-Z(\gamma)+\delta\tilde{N}(\alpha)|^5}g(\gamma)d\gamma,
		\end{aligned}
	\end{equation*}
	and
	\begin{equation*}
		\begin{aligned}
			I_{2,3}&=\int_{B_\eta}\!\!\frac{(h_n+\delta)^3\tilde{N}_1(\alpha)\tilde{N}_2(\alpha)\tilde{N}_3(\alpha)}{|Z(\alpha)-Z(\gamma)+(h_n+\delta)\tilde{N}(\alpha)|^5}g(\gamma)d\gamma -\int_{B_\eta}\frac{\delta^3\tilde{N}_1(\alpha)\tilde{N}_2(\alpha)\tilde{N}_3(\alpha)}{|Z(\alpha)-Z(\gamma)+\delta\tilde{N}(\alpha)|^5}g(\gamma)d\gamma.
		\end{aligned}
	\end{equation*}
	Next, we perform a further decomposition in $I_{2,1}$ to handle it:
	\begin{equation}\label{J3J7split}
	I_{2,1}=\sum_{k=3}^7J_k
	\end{equation}
	where
	\begin{equation*}
		\begin{aligned}
			J_3&=\int_{B_\eta}\!\!\frac{(Z_1(\alpha)-Z_1(\gamma)-\partial_\alpha Z_1(\alpha)\cdot(\alpha-\gamma))(Z_2(\alpha)-Z_2(\gamma))(h_n+\delta)\tilde{N}_3(\alpha)}{|Z(\alpha)-Z(\gamma)+(h_n+\delta)\tilde{N}(\alpha)|^5}g(\gamma)d\gamma \\
			&\quad-\int_{B_\eta}\frac{(Z_1(\alpha)-Z_1(\gamma)-\partial_\alpha Z_1(\alpha)\cdot(\alpha-\gamma))(Z_2(\alpha)-Z_2(\gamma))\delta\tilde{N}_3(\alpha)}{|Z(\alpha)-Z(\gamma)+\delta\tilde{N}(\alpha)|^5}g(\gamma)d\gamma,
		\end{aligned}
	\end{equation*}
	\begin{equation*}
		\begin{aligned}
			J_4&=\int_{B_\eta}\!\!\frac{\partial_\alpha Z_1(\alpha)\cdot(\alpha-\gamma)(Z_2(\alpha)-Z_2(\gamma)-\partial_\alpha Z_2(\alpha)\cdot(\alpha-\gamma))(h_n+\delta)\tilde{N}_3(\alpha)}{|Z(\alpha)-Z(\gamma)+(h_n+\delta)\tilde{N}(\alpha)|^5}g(\gamma)d\gamma \\
			&\quad-\int_{B_\eta}\frac{\partial_\alpha Z_1(\alpha)\cdot(\alpha-\gamma)(Z_2(\alpha)-Z_2(\gamma)-\partial_\alpha Z_2(\alpha)\cdot(\alpha-\gamma))\delta\tilde{N}_3(\alpha)}{|Z(\alpha)-Z(\gamma)+\delta\tilde{N}(\alpha)|^5}g(\gamma)d\gamma,
		\end{aligned}
	\end{equation*}
	\begin{equation*}
		\begin{aligned}
			J_5&=\int_{B_\eta}\!\!\frac{\partial_\alpha Z_1(\alpha)\cdot(\alpha-\gamma)\partial_\alpha Z_2(\alpha)\cdot(\alpha-\gamma)(h_n+\delta)\tilde{N}_3(\alpha)}{|Z(\alpha)-Z(\gamma)+(h_n+\delta)\tilde{N}(\alpha)|^5}(g(\gamma)-g(\alpha))d\gamma \\
			&\quad-\int_{B_\eta}\frac{\partial_\alpha Z_1(\alpha)\cdot(\alpha-\gamma)\partial_\alpha Z_2(\alpha)\cdot(\alpha-\gamma)\delta\tilde{N}_3(\alpha)}{|Z(\alpha)-Z(\gamma)+\delta\tilde{N}(\alpha)|^5}(g(\gamma)-g(\alpha))d\gamma,
		\end{aligned}
	\end{equation*}
	\begin{equation*}
		\begin{aligned}
			J_6&=(h_n+\delta)\tilde{N}_3(\alpha)g(\alpha)\int_{B_\eta}\partial_\alpha Z_1(\alpha)\cdot(\alpha-\gamma)\partial_\alpha Z_2(\alpha)\cdot(\alpha-\gamma)\big(\frac{1}{u_h^5}-\frac{1}{v_h^5}\big)d\gamma \\
			&\quad-\delta\tilde{N}_3(\alpha)g(\alpha)\int_{B_\eta}\partial_\alpha Z_1(\alpha)\cdot(\alpha-\gamma)\partial_\alpha Z_2(\alpha)\cdot(\alpha-\gamma)\big(\frac{1}{u_0^5}-\frac{1}{v_0^5}\big)d\gamma,
		\end{aligned}
	\end{equation*}
	(see \eqref{uhvh}) and
	\begin{equation*}
		\begin{aligned}
			J_7&=(h_n+\delta)\tilde{N}_3(\alpha)g(\alpha)\int_{B_\eta}\frac{\partial_\alpha Z_1(\alpha)\cdot(\alpha-\gamma)\partial_\alpha Z_2(\alpha)\cdot(\alpha-\gamma)}{|\partial_\alpha Z(\alpha)(\alpha-\gamma)+(h_n+\delta)\tilde{N}(\alpha)|^5}d\gamma \\
			&\quad-\delta\tilde{N}_3(\alpha)g(\alpha)\int_{B_\eta}\frac{\partial_\alpha Z_1(\alpha)\cdot(\alpha-\gamma)\partial_\alpha Z_2(\alpha)\cdot(\alpha-\gamma)}{|\partial_\alpha Z(\alpha)(\alpha-\gamma)+\delta\tilde{N}(\alpha)|^5}d\gamma.
		\end{aligned}
	\end{equation*}
	In the next step, a splitting gives
	$$
	J_{3}=J_{3,1}+J_{3,2}
	$$
	with
	\begin{equation*}
		\begin{aligned}
			J_{3,1}&=h_n\tilde{N}_3(\alpha)\int_{B_\eta}\!\!\frac{(Z_1(\alpha)-Z_1(\gamma)-\partial_\alpha Z_1(\alpha)\cdot(\alpha-\gamma))(Z_2(\alpha)-Z_2(\gamma))}{|Z(\alpha)-Z(\gamma)+(h_n+\delta)\tilde{N}(\alpha)|^5}g(\gamma)d\gamma 
		\end{aligned}
	\end{equation*}
	and
	\begin{equation*}
		\begin{aligned}
			J_{3,2}&=\delta\tilde{N}_3(\alpha)\int_{B_\eta}\!\!(Z_1(\alpha)-Z_1(\gamma)-\partial_\alpha Z_1(\alpha)\cdot(\alpha-\gamma))(Z_2(\alpha)-Z_2(\gamma))g(\gamma)\big(\frac{1}{u_h^5}-\frac{1}{u_0^5}\big)d\gamma.
		\end{aligned}
	\end{equation*}
	Bound \eqref{denominatorb2} provides the desired estimates as before:
	\begin{equation*}
		|J_{3,1}|\leq C  \|g\|_{L^\infty}\frac{\|\pa Z\|_{L^\infty}^2}{|\pa Z|_{\text{inf}}^5}\|\pa Z\|_{\dot{C}^\sigma}|h_n|^\sigma,
	\end{equation*}
	\begin{equation*}
		\begin{aligned}
			|J_{3,2}|&\leq C  \|g\|_{L^\infty}\frac{\|\pa Z\|_{L^\infty}^3}{|\pa Z|_{\text{inf}}^6}\|\pa Z\|_{\dot{C}^\sigma}\delta|h_n|\int_{B_{\eta}}\Big(\frac{|\alpha-\gamma|^{2+\sigma}(|\alpha-\gamma|^2+\delta^2)^{-\frac12}}{(|\alpha-\gamma|^2+(h_n+\delta)^2)^{\frac52}}\\
			&\qquad\qquad\qquad\qquad\qquad\qquad\qquad\qquad+\frac{|\alpha-\gamma|^{2+\sigma}(|\alpha-\gamma|^2+\delta^2)^{-\frac52}}{(|\alpha-\gamma|^2+(h_n+\delta)^2)^{\frac12}}\Big)d\gamma\\
			&\leq  C  \|g\|_{L^\infty}\frac{\|\pa Z\|_{L^\infty}^3}{|\pa Z|_{\text{inf}}^6}\|\pa Z\|_{\dot{C}^\sigma}|h_n|^{\sigma}.
		\end{aligned}
	\end{equation*}
	We are done with $J_3$. The terms $J_{4}$ and $J_{5}$ follow in a similar manner,
	\begin{equation*}
		|J_4|\leq  C  \|g\|_{L^\infty}\Big(\frac{\|\pa Z\|_{L^\infty}^2}{|\pa Z|_{\text{inf}}^5}+\frac{\|\pa Z\|_{L^\infty}^3}{|\pa Z|_{\text{inf}}^6}\Big)\|\pa Z\|_{\dot{C}^\sigma}|h_n|^{\sigma},
	\end{equation*}
	but for $J_5$ the bound is slightly different:
	\begin{equation*}
		|J_5|\leq  C  \| g\|_{\dot{C}^\sigma}\Big(\frac{\|\pa Z\|_{L^\infty}^3}{|\pa Z|_{\text{inf}}^5}+\frac{\|\pa Z\|_{L^\infty}^4}{|\pa Z|_{\text{inf}}^6}\Big)|h_n|^{\sigma}.
	\end{equation*}
	To handle $J_6$ we split it in two:
	\begin{equation*}
		J_6=J_{6,1}+J_{6,2},
	\end{equation*}
	where
	\begin{equation*}
		\begin{aligned}
			J_{6,1}&=h_n\tilde{N}_3(\alpha)g(\alpha)\int_{B_\eta}\partial_\alpha Z_1(\alpha)\cdot(\alpha-\gamma)\partial_\alpha Z_2(\alpha)\cdot(\alpha-\gamma)\big(\frac{1}{u_h^5}-\frac{1}{v_h^5}\big)d\gamma,
		\end{aligned}
	\end{equation*}
	and
	\begin{equation*}
		\begin{aligned}
			J_{6,2}&=\delta\tilde{N}_3(\alpha)g(\alpha)\int_{B_\eta}\partial_\alpha Z_1(\alpha)\cdot(\alpha-\gamma)\partial_\alpha Z_2(\alpha)\cdot(\alpha-\gamma)\Big(\frac{1}{u_h^5}-\frac{1}{u_0^5}-\big(\frac{1}{v_h^5}-\frac{1}{v_0^5}\big)\Big)d\gamma.
		\end{aligned}
	\end{equation*}
	As before it is possible to get
	\begin{equation*}
		|J_{6,1}|\leq C  \|g\|_{L^\infty}\frac{\|\pa Z\|_{L^\infty}^{11}}{|\pa Z|_{\text{inf}}^{14}}\|\pa Z\|_{\dot{C}^\sigma}|h_n|^{\sigma}.
	\end{equation*}
	The next term can be estimated similarly to  $I_{1,5}$ \eqref{I15} in order to obtain
	\begin{equation*}
		|J_{6,2}|\leq C  \|g\|_{L^\infty}\frac{\|\pa Z\|_{L^\infty}^{13}}{|\pa Z|_{\text{inf}}^{16}}\|\pa Z\|_{\dot{C}^\sigma}|h_n|^{\sigma},
	\end{equation*}
	and therefore the appropriate bound for $J_6$:
	\begin{equation*}
		|J_{6}|\leq C  \|g\|_{L^\infty}\frac{\|\pa Z\|_{L^\infty}^{13}}{|\pa Z|_{\text{inf}}^{16}}\|\pa Z\|_{\dot{C}^\sigma}|h_n|^{\sigma},
	\end{equation*}
We now show that the last term $J_7$ \eqref{J3J7split} is Lipschitz. The change of variables $\gamma\leftarrow (\alpha-\gamma)/(h_n+\delta)$ gives that
\begin{equation*}
	\begin{aligned}
		J_7&=\tilde{N}_3(\alpha)g(\alpha)\Big(\int_{|\gamma|\leq\frac{\eta}{h_n+\delta}}\frac{\partial_\alpha Z_1(\alpha)\cdot\gamma\partial_\alpha Z_2(\alpha)\cdot\gamma}{(|\partial_\alpha Z(\alpha)\gamma|^2+|\tilde{N}(\alpha)|^2)^{\frac52}}d\gamma-\int_{|\gamma|\leq\frac{\eta}{\delta}}\frac{\partial_\alpha Z_1(\alpha)\cdot\gamma\partial_\alpha Z_2(\alpha)\cdot\gamma}{(|\partial_\alpha Z(\alpha)\gamma|^2+|\tilde{N}(\alpha)|^2)^{\frac52}}d\gamma\Big).
	\end{aligned}
\end{equation*}
Define
\begin{equation*}
    \begin{aligned}
        F(z)&=\int_{|\gamma|\leq\frac{\eta}{z}}\frac{\partial_\alpha Z_1(\alpha)\cdot\gamma\partial_\alpha Z_2(\alpha)\cdot\gamma}{(|\partial_\alpha Z(\alpha)\gamma|^2+|\tilde{N}(\alpha)|^2)^{\frac52}}d\gamma,
    \end{aligned}
\end{equation*}
which, denoting $\hat{\gamma}=\frac{\gamma}{|\gamma|}$, can be written as follows
\begin{equation*}
    \begin{aligned}
        F(z)&=\int_{-\pi}^\pi\frac{\partial_\alpha Z_1(\alpha)\cdot\hat{\gamma}\partial_\alpha Z_2(\alpha)\cdot\hat{\gamma}}{|\tilde{N}(\alpha)|^5}\int_0^{\frac{\eta}{z}}\frac{r^3}{(|\tilde{N}(\alpha)|^{-2}|\partial_\alpha Z(\alpha)\hat{\gamma}|^2r^2+1)^{\frac52}}drd\hat{\gamma}. 
    \end{aligned}
\end{equation*}
If we denote further
\begin{equation}\label{J7aux}
    \begin{aligned}
        G(r,a)=\int_0^{\frac{1}{r}}\frac{\rho^3}{(a^2\rho^2+1)^{\frac52}}d\rho =-\frac{2+3a^2r^{-2}}{3a^4(1+a^2r^{-2})^{\frac32}},
    \end{aligned}
\end{equation}
we obtain that
\begin{equation*}
    \begin{aligned}
        |J_7|&\leq \|\pa Z\|_{L^\infty}\|g\|_{L^\infty}\Big|\int_{-\pi}^\pi\frac{\partial_\alpha Z_1(\alpha)\cdot\hat{\gamma}\partial_\alpha Z_2(\alpha)\cdot\hat{\gamma}}{|\tilde{N}(\alpha)|^5}\\
        &\qquad\times\Big(G(\frac{h_n+\delta}{\eta},|\tilde{N}(\alpha)|^{-1}|\partial_\alpha Z(\alpha)\hat{\gamma}|)-G(\frac{\delta}{\eta},|\tilde{N}(\alpha)|^{-1}|\partial_\alpha Z(\alpha)\hat{\gamma}|)\Big)d\hat{\gamma} \Big|.
    \end{aligned}
\end{equation*}
Since $|\frac{d}{dr}G(a,r)|\leq |a|^{-5}$, we obtain
\begin{equation*}
    \begin{aligned}
        |J_7|&\leq C\|g\|_{L^\infty}\frac{\|\pa Z\|_{L^\infty}^8}{|\pa Z|_{\inf}^{10}}\frac{|h_n|}{\eta}.
    \end{aligned}
\end{equation*}
	This yields the appropriate estimate for $I_{2,1}$.
	Next we handle $I_{2,2}$ with the splitting
	$$
	I_{2,2}=\sum_{k=8}^{11}J_{k},
	$$
	where
	\begin{equation*}
		\begin{aligned}
			J_{8}&=(h_n+\delta)^2\tilde{N}_2(\alpha)\tilde{N}_3(\alpha)\int_{B_\eta}\!\!\frac{(Z_1(\alpha)-Z_1(\gamma)-\pa Z_1(\alpha)(\alpha-\gamma))}{|Z(\alpha)-Z(\gamma)+(h_n+\delta)\tilde{N}(\alpha)|^5}g(\gamma)d\gamma \\
			&\quad-\delta^2\tilde{N}_2(\alpha)\tilde{N}_3(\alpha)\int_{B_\eta}\frac{(Z_1(\alpha)-Z_1(\gamma)-\pa Z_1(\alpha)(\alpha-\gamma))}{|Z(\alpha)-Z(\gamma)+\delta\tilde{N}(\alpha)|^5}g(\gamma)d\gamma,
		\end{aligned}
	\end{equation*}
	\begin{equation*}
		\begin{aligned}
			J_{9}&=(h_n+\delta)^2\tilde{N}_2(\alpha)\tilde{N}_3(\alpha)\int_{B_\eta}\!\!\frac{\pa Z_1(\alpha)(\alpha-\gamma)}{|Z(\alpha)-Z(\gamma)+(h_n+\delta)\tilde{N}(\alpha)|^5}(g(\gamma)-g(\alpha))d\gamma \\
			&\quad-\delta^2\tilde{N}_2(\alpha)\tilde{N}_3(\alpha)\int_{B_\eta}\frac{\pa Z_1(\alpha)(\alpha-\gamma)}{|Z(\alpha)-Z(\gamma)+\delta\tilde{N}(\alpha)|^5}(g(\gamma)-g(\alpha))d\gamma,
		\end{aligned}
	\end{equation*}
	\begin{equation*}
		\begin{aligned}
			J_{10}&=(h_n+\delta)^2\tilde{N}_2(\alpha)\tilde{N}_3(\alpha)g(\alpha)\int_{B_\eta}\!\!\pa Z_1(\alpha)(\alpha-\gamma)\big(\frac{1}{u_h^5}-\frac{1}{v^5_h}\big)d\gamma \\
			&\quad-\delta^2\tilde{N}_2(\alpha)\tilde{N}_3(\alpha)g(\alpha)\int_{B_\eta}\pa Z_1(\alpha)(\alpha-\gamma)\big(\frac{1}{u_0^5}-\frac{1}{v^5_0}\big)d\gamma,
		\end{aligned}
	\end{equation*}
	(see \eqref{uhvh}) and
	\begin{equation*}
		\begin{aligned}
			J_{11}&=(h_n+\delta)^2\tilde{N}_2(\alpha)\tilde{N}_3(\alpha)g(\alpha)\int_{B_\eta}\!\!\frac{\pa Z_1(\alpha)(\alpha-\gamma)}{|\pa Z(\alpha)(\alpha-\gamma)+(h_n+\delta)\tilde{N}(\alpha)|^5}d\gamma \\
			&\quad-\delta^2\tilde{N}_2(\alpha)\tilde{N}_3(\alpha)g(\alpha)\int_{B_\eta}\frac{\pa Z_1(\alpha)(\alpha-\gamma)}{|\pa Z(\alpha)(\alpha-\gamma)+\delta\tilde{N}(\alpha)|^5}d\gamma.
		\end{aligned}
	\end{equation*}
	The term $J_{8}$ can be decomposed further to get
	$$
	J_8=J_{8,1}+J_{8,2},
	$$
	where
	\begin{equation*}
		\begin{aligned}
			J_{8,1}&=(h_n^2+2h_n\delta)\tilde{N}_2(\alpha)\tilde{N}_3(\alpha)\int_{B_\eta}\!\!\frac{Z_1(\alpha)-Z_1(\gamma)-\pa Z_1(\alpha)(\alpha-\gamma)}{|Z(\alpha)-Z(\gamma)+(h_n+\delta)\tilde{N}(\alpha)|^5}g(\gamma)d\gamma
		\end{aligned}
	\end{equation*}    
	and
	\begin{equation*}
		\begin{aligned}
			J_{8,2}&=\delta^2\tilde{N}_2(\alpha)\tilde{N}_3(\alpha)\int_{B_\eta}(Z_1(\alpha)-Z_1(\gamma)-\pa Z_1(\alpha)(\alpha-\gamma))g(\gamma)\big(\frac1{u_h^5}-\frac1{u_0^5}\big)d\gamma.
		\end{aligned}
	\end{equation*}
Then, it is possible to bound as follows
\begin{equation*}
	|J_{8,1}|\leq C \|g\|_{L^\infty}\frac{\|\pa Z\|_{L^\infty}^2}{|\pa Z|_{\text{inf}}^5}\|\pa Z\|_{\dot{C}^\sigma}\frac{|h_n|(2|h_n+\delta|+|h_n|)}{|h_n+\delta|^{2-\sigma}}\leq C \|g\|_{L^\infty}\frac{\|\pa Z\|_{L^\infty}^2}{|\pa Z|_{\text{inf}}^5}\|\pa Z\|_{\dot{C}^\sigma}|h_n|^{\sigma},
\end{equation*}
to get the desired estimate. Similarly as before, the next term is approached:
\begin{equation*}
	\begin{aligned}
		|J_{8,2}|&\leq C  \|g\|_{L^\infty}\frac{\|\pa Z\|_{L^\infty}^8}{|\pa Z|_{\text{inf}}^{11}}\|\pa Z\|_{\dot{C}^\sigma}\delta^2|h_n|\int_{B_{\eta}}\Big(\frac{|\alpha-\gamma|^{1+\sigma}(|\alpha-\gamma|^2+\delta^2)^{-\frac12}}{(|\alpha-\gamma|^2+(h_n+\delta)^2)^{\frac52}}\\
		&\qquad\qquad\qquad\qquad\qquad\qquad\qquad\qquad+\frac{|\alpha-\gamma|^{1+\sigma}(|\alpha-\gamma|^2+\delta^2)^{-\frac52}}{(|\alpha-\gamma|^2+(h_n+\delta)^2)^{\frac12}}\Big)d\gamma\\
		&\leq  C  \|g\|_{L^\infty}\frac{\|\pa Z\|_{L^\infty}^8}{|\pa Z|_{\text{inf}}^{11}}\|\pa Z\|_{\dot{C}^\sigma}\Big(\frac{\delta|h_n|}{|h_n\!+\!\delta|^{2-\sigma}}\!+\!\frac{|h_n|}{|h_n\!+\!\delta|^{1-\sigma}}\Big)\leq
		C  \|g\|_{L^\infty}\frac{\|\pa Z\|_{L^\infty}^8}{|\pa Z|_{\text{inf}}^{11}}\|\pa Z\|_{\dot{C}^\sigma}|h_n|^{\sigma}.
	\end{aligned}
\end{equation*}
It yields the desired estimate for $J_8$:
\begin{equation*}
	|J_{8}|\leq C \|g\|_{L^\infty}\frac{\|\pa Z\|_{L^\infty}^8}{|\pa Z|_{\text{inf}}^{11}}\|\pa Z\|_{\dot{C}^\sigma}|h_n|^{\sigma}.
\end{equation*}
The next term can be handled analogously, obtaining the bound below
\begin{equation*}
	|J_9|\leq  C  \| g\|_{\dot{C}^\sigma}\frac{\|\pa Z\|_{L^\infty}^9}{|\pa Z|_{\text{inf}}^{11}}|h_n|^{\sigma}.
\end{equation*}
We continue dealing with $J_{10}$, which can be decomposed as before to obtain
\begin{equation*}
	|J_{10}|\leq C  \|g\|_{L^\infty}\frac{\|\pa Z\|_{L^\infty}^{10}}{|\pa Z|_{\text{inf}}^{13}}\|\pa Z\|_{\dot{C}^\sigma}|h_n|^{\sigma}.
\end{equation*}
The last term in the splitting is prepared to be integrated explicitly. In this case it is easy to check that it is zero: 
$$
J_{11}=0.
$$	
Gathering the last four estimates provides
the appropriate estimate for $I_{2,2}$. It remains to handle $I_{2,3}$. A further decomposition yields 	
$$
I_{2,3}=J_{12}+J_{13}+J_{14},
$$
where
\begin{equation*}
	\begin{aligned}
		J_{12}=&(h_n+\delta)^3\tilde{N}_1(\alpha)\tilde{N}_2(\alpha)\tilde{N}_3(\alpha)\int_{B_\eta}\!\!\frac{(g(\gamma)-g(\alpha))d\gamma}{|Z(\alpha)-Z(\gamma)+(h_n+\delta)\tilde{N}(\alpha)|^5}\\ &-\delta^3\tilde{N}_1(\alpha)\tilde{N}_2(\alpha)\tilde{N}_3(\alpha)\int_{B_\eta}\frac{(g(\gamma)-g(\alpha))d\gamma}{|Z(\alpha)-Z(\gamma)+\delta\tilde{N}(\alpha)|^5},
	\end{aligned}
\end{equation*}
\begin{equation*}
	\begin{aligned}
		J_{13}=&(h_n+\delta)^3\tilde{N}_1(\alpha)\tilde{N}_2(\alpha)\tilde{N}_3(\alpha)g(\alpha)\int_{B_\eta}\!\!\big(\frac{1}{u_h^5}-\frac{1}{v_h^5}\big)d\gamma\\ &-\delta^3\tilde{N}_1(\alpha)\tilde{N}_2(\alpha)\tilde{N}_3(\alpha)g(\alpha)\int_{B_\eta}\!\!\big(\frac{1}{u_0^5}-\frac{1}{v_0^5}\big)d\gamma,
	\end{aligned}
\end{equation*}
and
	\begin{equation*}
	\begin{aligned}
		J_{14}&=(h_n+\delta)^3\tilde{N}_1(\alpha)\tilde{N}_2(\alpha)\tilde{N}_3(\alpha)g(\alpha)\int_{B_\eta}\!\!\frac{d\gamma}{(|\pa Z(\alpha)(\alpha-\gamma)|^2+(h_n+\delta)^2|\tilde{N}(\alpha)|^2)^{\frac52}} \\
		&\quad-\delta^3\tilde{N}_1(\alpha)\tilde{N}_2(\alpha)\tilde{N}_3(\alpha)g(\alpha)\int_{B_\eta}\frac{d\gamma}{(|\pa Z(\alpha)(\alpha-\gamma)|^2+\delta^2|\tilde{N}(\alpha)|^2)^{\frac52}}.
	\end{aligned}
\end{equation*}
A further decomposition helps to deal with $J_{12}$:
$$
J_{12}=J_{12,1}+J_{12,2},
$$
where
\begin{equation*}
	\begin{aligned}
		J_{12,1}=&(h_n^3+3h_n\delta(h_n+\delta))\tilde{N}_1(\alpha)\tilde{N}_2(\alpha)\tilde{N}_3(\alpha)\int_{B_\eta}\!\!\frac{(g(\gamma)-g(\alpha))d\gamma}{|Z(\alpha)-Z(\gamma)+(h_n+\delta)\tilde{N}(\alpha)|^5},
	\end{aligned}
\end{equation*}
	and
\begin{equation*}
	\begin{aligned}
		J_{12,2}=&\delta^3\tilde{N}_1(\alpha)\tilde{N}_2(\alpha)\tilde{N}_3(\alpha)\int_{B_\eta}\!\!
		(g(\gamma)-g(\alpha))\big(\frac{1}{u_h^5}-\frac{1}{u_0^5}\big)d\gamma.
	\end{aligned}
\end{equation*}	
Next, it is possible to bound as follows
\begin{equation*}
	|J_{12,1}|\leq C \frac{\|\pa Z\|_{L^\infty}^3}{|\pa Z|_{\text{inf}}^5}\|g\|_{\dot{C}^\sigma}\frac{|h_n|(|h_n|^2\!+\!3|h_n\!+\!\delta|^2\!+\!3|h_n||h_n\!+\!\delta|)}{|h_n+\delta|^{3-\sigma}}\leq C \frac{\|\pa Z\|_{L^\infty}^3}{|\pa Z|_{\text{inf}}^5}\|g\|_{\dot{C}^\sigma}|h_n|^{\sigma},
\end{equation*}
\begin{equation*}
	\begin{aligned}
		|J_{12,2}|&\leq C  \frac{\|\pa Z\|_{L^\infty}^4}{|\pa Z|_{\text{inf}}^6}\|g\|_{\dot{C}^\sigma}\delta^3|h_n|\int_{B_{\eta}}|\alpha-\gamma|^{\sigma}\Big(\frac{(|\alpha-\gamma|^2+\delta^2)^{-\frac12}}{(|\alpha-\gamma|^2+(h_n+\delta)^2)^{\frac52}}\\
		&\qquad\qquad\qquad\qquad\qquad\qquad\qquad\qquad\quad+\frac{(|\alpha-\gamma|^2+\delta^2)^{-\frac52}}{(|\alpha-\gamma|^2+(h_n+\delta)^2)^{\frac12}}\Big)d\gamma\\
		&\leq  C  \frac{\|\pa Z\|_{L^\infty}^4}{|\pa Z|_{\text{inf}}^6}\|g\|_{\dot{C}^\sigma}\Big(\frac{\delta^2|h_n|}{|h_n\!+\!\delta|^{3-\sigma}}\!+\!\frac{|h_n|}{|h_n\!+\!\delta|^{1-\sigma}}\Big)\leq
		C \frac{\|\pa Z\|_{L^\infty}^4}{|\pa Z|_{\text{inf}}^6}\|g\|_{\dot{C}^\sigma}|h_n|^{\sigma}.
	\end{aligned}
\end{equation*}
and to get the appropriate estimate for $J_{12}$. A similar splitting for $J_{13}$ allows to obtain
\begin{equation*}
	|J_{13}|\leq C  \|g\|_{L^\infty}\Big(\frac{\|\pa Z\|_{L^\infty}^4}{|\pa Z|_{\text{inf}}^6}+\frac{\|\pa Z\|_{L^\infty}^5}{|\pa Z|_{\text{inf}}^7}+\frac{\|\pa Z\|_{L^\infty}^6}{|\pa Z|_{\text{inf}}^8}\Big)\|\pa Z\|_{\dot{C}^\sigma}|h_n|^{\sigma}.
\end{equation*}
Similarly as we did for $J_7$ in \eqref{J7aux}, after a change of variables, the radial part of the integrals in $J_{14}$ can be integrated to obtain the desired bound. It yields the appropriate estimate for $I_{2,3}$. We therefore complete controlling the term $I_2$ \eqref{I2split},
	\begin{equation}\label{I2bound}
	\begin{aligned}
	    |I_2|&\leq C \minspace \|g\|_{L^\infty}\|\pa Z\|_{L^\infty}^2\|F(Z)\|_{L^\infty}^5(1+\|\pa Z\|_{L^\infty}^{8}\|F(Z)\|_{L^\infty}^{8})\|\pa Z\|_{\dot{C}^\sigma}|h_n|^\sigma\\
	        &\quad+C \minspace \|g\|_{\dot{C}^\sigma}\|\pa Z\|_{L^\infty}^3\|F(Z)\|_{L^\infty}^6(1+\|\pa Z\|_{L^\infty}^6\|F(Z)\|_{L^\infty})^6|h_n|^\sigma.
	\end{aligned}
	\end{equation}
	Together with \eqref{I1bound}, it provides the desired bound for $I$ \eqref{Isplit}. Recalling \eqref{gbound}, this completes the proof of the Case 1,
	\begin{equation}\label{normal_bound}
	    \begin{aligned}
	        |S(f)(Z(\alpha)+(\delta+h_n)\tilde{N}(\alpha))-S(f)(Z(\alpha)+\delta &\tilde{N}(\alpha))|\leq\\ 
	        P(\|\pa Z\|_{L^\infty}\!&+\!\|F(Z)\|_{L^\infty})\big(\|f\|_{L^\infty}\|\pa Z\|_{\dot{C}^\sigma}\!+\! \|f\|_{\dot{C}^\sigma}\big)|h_n|^\sigma.
	    \end{aligned}
	\end{equation}

	\vspace{.2cm}

	\noindent\textbf{Case 2: Tangential direction.}
	Here we consider the first difference in \eqref{normaltangentSplit}. We recall that in the case we are dealing with, $|h_\tau|\leq \frac14\frac{|\pa Z|_{\inf}}{\|\pa Z\|_{L^\infty}}\delta$, we have that $\delta+h_n>0$. Hence, for simplicity in notation we do the estimate for some $\delta>0$, and we will later apply it with $\delta+h_n$. We keep the notation $x=Z(\alpha)+\delta \tilde{N}(\alpha)$. First, we split as before
	\begin{equation}\label{holder_tan}
		\begin{aligned}
			S(f)(x+h_\tau\cdot\pa Z(\alpha))-S(f)(x)&=\pv\int_{Z(B_\eta)}\Big(k(x+h_\tau\cdot\pa Z(\alpha)-y)-k(x-y)\Big)f(y)dS(y)\\
			&+\int_{\partial D \smallsetminus Z(B_\eta)}\Big(k(x+h_\tau\cdot\pa Z(\alpha)-y)-k(x-y)\Big)f(y)dS(y)\\
			&=I+II,
		\end{aligned}
	\end{equation}
	where 
	\begin{equation*}
		B_\eta=\{\gamma\in\mathbb{R}^2:|\alpha-\gamma|<\eta\}.
	\end{equation*}
	The second term is again more regular, 
	\begin{equation*}
		\begin{aligned}
			|II|\leq \frac{C}{\eta^3|\pa Z|_{\inf}}|\partial D|\|f\|_{L^\infty}|h|. 
		\end{aligned}
	\end{equation*}
	The first term is given by
	\begin{equation*}
		\begin{aligned}
			I&=\int_{B_\eta} \frac{\prod_{j=1}^3(x_j+h_\tau\cdot\pa Z_j(\alpha)-Z_j(\gamma))}{|x+h_\tau\cdot\pa Z(\alpha)-Z(\gamma)|^5}g(\gamma)d\gamma-\int_{B_\eta} \frac{\prod_{j=1}^3(x_j-Z_j(\gamma))}{|x-Z(\gamma)|^5}g(\gamma)d\gamma,
		\end{aligned}    
	\end{equation*}
	and we decompose it as follows
	\begin{equation}\label{Isplit_tangent}
		I=\sum_{i=1}^6 I_{i},
	\end{equation}
	where
	\begin{equation}\label{I1_tangent}
		\begin{aligned}
			I_1&= \int_{B_\eta}\!\!\! \frac{\prod_{j=1}^2(Z_j(\alpha)\!-\!Z_j(\gamma)\!+\!h_\tau\cdot\pa Z_j(\alpha)\!+\!\delta\tilde{N}_j(\alpha))(Z_3(\alpha)\!-\!Z_3(\gamma)\!-\!(\alpha\!-\!\gamma)\cdot \pa Z_3(\alpha))}{|Z(\alpha)-Z(\gamma)+h_\tau\cdot\pa Z(\alpha)+\delta \tilde{N}(\alpha)|^5}g(\gamma)d\gamma\\
			&\quad-\int_{B_\eta} \frac{\prod_{j=1}^2(Z_j(\alpha)-Z_j(\gamma)+\delta\tilde{N}_j(\alpha))(Z_3(\alpha)-Z_3(\gamma)-(\alpha-\gamma)\cdot \pa Z_3(\alpha))}{|Z(\alpha)-Z(\gamma)+\delta \tilde{N}(\alpha)|^5}g(\gamma)d\gamma,
		\end{aligned}
	\end{equation}
	\begin{equation}\label{I2_tangent}
		\begin{aligned}
			I_2&= \int_{B_\eta}\!\!\! \frac{((\alpha-\gamma+h_\tau)\cdot\pa Z_3(\alpha)+\delta\tilde{N}_3(\alpha))(Z_2(\alpha)-Z_2(\gamma)-(\alpha-\gamma)\cdot\pa Z_2(\alpha))}{|Z(\alpha)-Z(\gamma)+h_\tau\cdot\pa Z(\alpha)+\delta \tilde{N}(\alpha)|^5}\\
			&\qquad\times(Z_1(\alpha)-Z_1(\gamma)+h_\tau\cdot\pa Z_1(\alpha)+\delta\tilde{N}_1(\alpha))g(\gamma)d\gamma\\
			& \quad-\int_{B_\eta} \frac{((\alpha-\gamma)\cdot\pa Z_3(\alpha)+\delta\tilde{N}_3(\alpha))(Z_2(\alpha)-Z_2(\gamma)-(\alpha-\gamma)\cdot\pa Z_2(\alpha))}{|Z(\alpha)-Z(\gamma)+\delta \tilde{N}(\alpha)|^5}\\
			&\qquad\times (Z_1(\alpha)-Z_1(\gamma)+\delta\tilde{N}_1(\alpha)) g(\gamma)d\gamma,
		\end{aligned}
	\end{equation}
	\begin{equation}\label{I3_tangent}
		\begin{aligned}
			I_3&= \int_{B_\eta}\!\!\! \frac{\prod_{j=2}^3((\alpha-\gamma+h_\tau)\cdot \pa Z_j(\alpha)+\delta\tilde{N}_j(\alpha))(Z_1(\alpha)-Z_1(\gamma)-(\alpha-\gamma)\cdot\pa Z_1(\alpha))}{|Z(\alpha)-Z(\gamma)+h_\tau\cdot\pa Z(\alpha)+\delta \tilde{N}(\alpha)|^5}g(\gamma)d\gamma\\
			&\quad-\int_{B_\eta} \frac{\prod_{j=2}^3((\alpha-\gamma)\cdot \pa Z_j(\alpha)+\delta\tilde{N}_j(\alpha))(Z_1(\alpha)-Z_1(\gamma)-(\alpha-\gamma)\cdot\pa Z_1(\alpha))}{|Z(\alpha)-Z(\gamma)+\delta \tilde{N}(\alpha)|^5}g(\gamma)d\gamma,
		\end{aligned}
	\end{equation}
	\begin{equation}\label{I4_tangent}
		\begin{aligned}
			I_4&= \int_{B_\eta}\!\!\! \frac{\prod_{j=1}^3((\alpha-\gamma+h_\tau)\cdot \pa Z_j(\alpha)+\delta\tilde{N}_j(\alpha))(g(\gamma)-g(\alpha))}{|Z(\alpha)-Z(\gamma)+h_\tau\cdot\pa Z(\alpha)+\delta \tilde{N}(\alpha)|^5}d\gamma\\
			&\quad-\int_{B_\eta} \frac{\prod_{j=1}^3((\alpha-\gamma)\cdot \pa Z_j(\alpha)+\delta\tilde{N}_j(\alpha))(g(\gamma)-g(\alpha))}{|Z(\alpha)-Z(\gamma)+\delta \tilde{N}(\alpha)|^5}d\gamma,
		\end{aligned}
	\end{equation}
	\begin{equation}\label{I5_tangent}
		\begin{aligned}
			I_5&= \int_{B_\eta} g(\alpha)\prod_{j=1}^3((\alpha-\gamma+h_\tau)\cdot \pa Z_j(\alpha)+\delta\tilde{N}_j(\alpha))\\
			&\qquad\times\Big(\frac{1}{|Z(\alpha)-Z(\gamma)+h_\tau\cdot\pa Z(\alpha)+\delta \tilde{N}(\alpha)|^5}-\frac{1}{\big(|\pa Z(\alpha)(\alpha-\gamma+h_\tau)|^2+\delta^2|\tilde{N}(\alpha)|^2\big)^{\frac52}}\Big)d\gamma\\
			&\quad-\int_{B_\eta} g(\alpha)\prod_{j=1}^3((\alpha-\gamma)\cdot \pa Z_j(\alpha)+\delta\tilde{N}_j(\alpha))\\
			&\quad\qquad\times \Big(\frac{1}{|Z(\alpha)-Z(\gamma)+\delta \tilde{N}(\alpha)|^5}-\frac{1}{\big(|\pa Z(\alpha)(\alpha-\gamma)|^2+\delta^2\tilde{N}(\alpha)|^2\big)^{\frac52}}\Big)d\gamma,
		\end{aligned}
	\end{equation}
	\begin{equation}\label{I6_tangent}
		\begin{aligned}
			I_6&= g(\alpha)\int_{B_\eta} \frac{\prod_{j=1}^3((\alpha-\gamma+h_\tau)\cdot \pa Z_j(\alpha)+\delta\tilde{N}_j(\alpha))}{\big(|\pa Z(\alpha)(\alpha-\gamma+h_\tau)|^2+\delta^2|\tilde{N}(\alpha)|^2\big)^{\frac52}}d\gamma\\
			&\quad-g(\alpha)\int_{B_\eta} \frac{\prod_{j=1}^3((\alpha-\gamma)\cdot \pa Z_j(\alpha)+\delta\tilde{N}_j(\alpha))}{\big(|\pa Z(\alpha)(\alpha-\gamma)|^2+\delta^2|\tilde{N}(\alpha)|^2\big)^{\frac52}}d\gamma.
		\end{aligned}
	\end{equation}
	We proceed to estimate each of these terms. 
	We first estimate the denominator as follows
	\begin{equation*}
		\begin{aligned}
			D&=|x-Z(\gamma)+h_{\tau}\cdot \pa Z(\alpha)|^2\\
			&=|Z(\alpha)-Z(\gamma)|^2+|h_\tau\cdot\pa Z(\alpha)|^2+\delta^2|\tilde{N}(\alpha)|^2\\
			&\quad+2h_\tau \cdot\pa Z(\alpha)\cdot(Z(\alpha)-Z(\gamma))+2\delta(Z(\alpha)-Z(\gamma)-\pa Z(\alpha)(\alpha-\gamma))\cdot\tilde{N}(\alpha)\\
			&\geq \frac{|Z(\alpha)-Z(\gamma)|^2}{|\alpha-\gamma|^2} |\alpha-\gamma|^2+|h_\tau\pa Z(\alpha)|^2+\delta^2|\tilde{N}(\alpha)|^2\\
			&\quad-2\delta|\alpha-\gamma|^{1+\sigma}\|\pa Z\|_{\dot{C}^\sigma}|\tilde{N}(\alpha)|-2|\alpha-\gamma||h_\tau\cdot\pa Z(\alpha)|\frac{|Z(\alpha)-Z(\gamma)|}{|\alpha-\gamma|}.
		\end{aligned}
	\end{equation*}
	The last two terms satisfy that
	\begin{equation*}
		\begin{aligned}
			2\delta|\alpha-\gamma|^{1+\sigma}\|\pa Z\|_{\dot{C}^\sigma}|\tilde{N}(\alpha)|&\leq\frac{1-\sigma}2\delta^2|\tilde{N}(\alpha)|^2\\
			&\hspace{-0.3cm}+ \frac{1+\sigma}2 2^{\frac{2}{1+\sigma}}\delta^{\frac{2\sigma}{1+\sigma}}\|\tilde{N}\|_{L^\infty}^{\frac{2\sigma}{1+\sigma}}\|\pa Z\|_{\dot{C}}^{\frac{2}{1+\sigma}}|\alpha-\gamma|^2,
		\end{aligned}
	\end{equation*}
	and
	\begin{equation*}
		2|\alpha-\gamma||h_\tau\cdot\pa Z(\alpha)|\frac{|Z(\alpha)-Z(\gamma)|}{|\alpha-\gamma|}\leq \frac14|\alpha-\gamma|^2\frac{|Z(\alpha)-Z(\gamma)|^2}{|\alpha-\gamma|^2}+4|h_\tau\cdot\pa Z(\alpha)|^2.
	\end{equation*}
	The choice of the cutoff for $\delta$ \eqref{delta_cutoff} and the fact that we are in the case where $|h_\tau|\leq \frac14\frac{|\pa Z|_{\inf}}{\|\pa Z\|_{L^\infty}}\delta$
	provides that
	\begin{equation}\label{denominatorb2_tan}
		\begin{aligned}
			D&\geq \frac12 |\pa Z|_{\text{inf}}^2\Big( |\alpha-\gamma|^2+\frac12\delta^2\Big).
		\end{aligned}
	\end{equation}
	Next, we split $I_1$ \eqref{Isplit_tangent} as follows
	\begin{equation}\label{I1split_tangent}
		\begin{aligned}
			I_1=I_{1,1}+I_{1,2}+I_{1,3},
		\end{aligned}
	\end{equation}
	with
	\begin{equation*}
		\begin{aligned}
			I_{1,1}&=h_\tau\cdot\pa Z_1(\alpha)\\
			&\quad\times\int_{B_\eta}\!\! \frac{(Z_2(\alpha)\!-\!Z_2(\gamma)+h_\tau\cdot\pa Z_2(\alpha)\!+\!\delta\tilde{N}_2(\alpha))(Z_3(\alpha)\!-\!Z_3(\gamma)\!-\!(\alpha\!-\!\gamma)\cdot\pa Z_3(\alpha))}{|Z(\alpha)-Z(\gamma)+h_\tau\cdot\pa Z(\alpha)+\delta \tilde{N}(\alpha)|^5}g(\gamma)d\gamma,
		\end{aligned}
	\end{equation*}
	\begin{equation*}
		\begin{aligned}
			I_{1,2}&\!=\!h_\tau\cdot\pa Z_2(\alpha)\!\int_{B_\eta}\!\! \frac{(Z_1(\alpha)\!-\!Z_1(\gamma)\!+\!\delta\tilde{N}_1(\alpha))(Z_3(\alpha)\!-\!Z_3(\gamma)\!-\!(\alpha\!-\!\gamma)\cdot\pa Z_3(\alpha))}{|Z(\alpha)-Z(\gamma)+h_\tau\cdot\pa Z(\alpha)+\delta \tilde{N}(\alpha)|^5}g(\gamma)d\gamma,
		\end{aligned}
	\end{equation*}
	\begin{equation*}
		\begin{aligned}
			I_{1,3}&=\int_{B_\eta} \prod_{j=1}^2(Z_j(\alpha)\!-\!Z_j(\gamma)\!+\!\delta\tilde{N}_j(\alpha))(Z_3(\alpha)\!-\!Z_3(\gamma)\!-\!(\alpha\!-\!\gamma)\cdot\pa Z_3(\alpha))g(\gamma)\\
			&\qquad\Big(\frac{1}{|Z(\alpha)-Z(\gamma)+h_\tau\cdot\pa Z(\alpha)+\delta \tilde{N}(\alpha)|^5}-\frac{1}{|Z(\alpha)-Z(\gamma)+\delta \tilde{N}(\alpha)|^5}\Big)d\gamma.
		\end{aligned}
	\end{equation*}
	We have that
	\begin{equation*}
		\begin{aligned}
			|I_{1,1}|&\leq |h_\tau|\frac{|\pa Z_1(\alpha)|\|\pa Z_2\|_{L^\infty}\|\pa Z_3\|_{\dot{C}^\sigma}\|g\|_{L^\infty}}{|\pa Z|_{\inf}^5}
			\int_{B_\eta} \frac{(|\alpha-\gamma|+|h_\tau|+\delta)|\alpha-\gamma|^{1+\sigma}}{\Big(|\alpha-\gamma|^2+\frac12\delta^2\Big)^5}d\gamma\\
			&\leq\frac{|h_\tau|}{\delta^{1-\sigma}}\frac{\|\pa Z\|_{L^\infty}^2\|\pa Z\|_{\dot{C}^\sigma}\|g\|_{L^\infty}}{|\pa Z|_{\inf}^5}
			\int_{\mathbb{R}^2} \frac{(|\gamma|+\frac{|h_\tau|}{\delta}+1)|\gamma|^{1+\sigma}}{\Big(|\gamma|^2+\frac12\Big)^5}d\gamma,
		\end{aligned}
	\end{equation*}
	and hence, recalling that we are dealing with the case where \eqref{htau_menor_delta} holds, we obtain that    
	\begin{equation*}
		\begin{aligned}
			|I_{1,1}|\leq  C\frac{\|\pa Z\|_{L^\infty}^{1+\sigma}\|g\|_{L^\infty}}{|\pa Z|_{\inf}^{4+\sigma}}\|\pa Z\|_{\dot{C}^\sigma}|h_\tau|^\sigma.
		\end{aligned}
	\end{equation*}
	Analogously, we find that
	\begin{equation*}
		\begin{aligned}
			|I_{1,2}|\leq  C\frac{\|\pa Z\|_{L^\infty}^{1+\sigma}\|g\|_{L^\infty}}{|\pa Z|_{\inf}^{4+\sigma}}\|\pa Z\|_{\dot{C}^\sigma}|h_\tau|^\sigma.
		\end{aligned}
	\end{equation*}
	If we denote
	\begin{equation}\label{uh_tan}
		\begin{aligned}
			u_h&=|Z(\alpha)-Z(\gamma)+h_\tau\cdot\pa Z(\alpha)+\delta\tilde{N}(\alpha)|,\\
			u_0&=|Z(\alpha)-Z(\gamma)+\delta\tilde{N}(\alpha)|,
		\end{aligned}
	\end{equation}
	we have that
	\begin{equation*}
		\begin{aligned}
			\frac{1}{u_h^5}-\frac{1}{u_0^5}=G(u_h,u_0)(u_0^2-u_h^2),
		\end{aligned}
	\end{equation*}
	where $G$ is defined in \eqref{Guhvh}. Since
	\begin{equation}\label{uhu0_tan}
	\begin{aligned}
    u_0^2-u_h^2&=-|h_\tau\cdot\pa Z(\alpha)|^2-2(Z(\alpha)-Z(\gamma)+\delta \tilde{N}(\alpha))\cdot(h_\tau\cdot\pa Z(\alpha))\\
	&=-|h_\tau\cdot\pa Z(\alpha)|^2-2(Z(\alpha)-Z(\gamma))\cdot h_\tau\cdot\pa Z(\alpha),
	\end{aligned}
\end{equation}
	we obtain that
	\begin{equation*}
		\begin{aligned}
			|I_{1,3}|&\leq C\frac{\|\pa Z\|_{L^\infty}^{3+\sigma}\|g\|_{L^\infty}}{|\pa Z|_{\inf}^{6+\sigma}}\|\pa Z\|_{\dot{C}^\sigma}|h_\tau|^\sigma,
		\end{aligned}
	\end{equation*}
	and thus
	\begin{equation}\label{I1bound_tan}
		|I_1|\leq C\frac{\|\pa Z\|_{L^\infty}^{1+\sigma}\|g\|_{L^\infty}}{|\pa Z|_{\inf}^{4+\sigma}}\Big(1+\frac{\|\pa Z\|_{L^\infty}^2}{|\pa Z|_{\inf}^{2}}\Big)\|\pa Z\|_{\dot{C}^\sigma}|h_\tau|^\sigma.
	\end{equation}
	The terms $I_2$ and $I_3$ \eqref{I1split_tangent} follow similarly and share the same bound
	\begin{equation}\label{I2I3bound_tan}
		|I_2|+|I_3|\leq C\frac{\|\pa Z\|_{L^\infty}^{1+\sigma}\|g\|_{L^\infty}}{|\pa Z|_{\inf}^{4+\sigma}}\Big(1+\frac{\|\pa Z\|_{L^\infty}^2}{|\pa Z|_{\inf}^{2}}\Big)\|\pa Z\|_{\dot{C}^\sigma}|h_\tau|^\sigma,
	\end{equation}
	while for $I_4$ we find that
	\begin{equation}\label{I4bound_tan}
		|I_4|\leq C\frac{\|\pa Z\|_{L^\infty}^{2+\sigma}}{|\pa Z|_{\inf}^{4+\sigma}}\Big(1+\frac{\|\pa Z\|_{L^\infty}^2}{|\pa Z|_{\inf}^{2}}\Big)\|g\|_{\dot{C}^\sigma}|h_\tau|^\sigma.
	\end{equation}
	We proceed to estimate $I_{5}$. We recall the notation \eqref{uh_tan} and define
	\begin{equation}\label{vh_tan}
		\begin{aligned}
			v_h&=\big(|\pa Z(\alpha)(\alpha-\gamma+h_\tau)|^2+\delta^2|\tilde{N}(\alpha)|^2\big)^{\frac12},\\
			v_0&=\big(|\pa Z(\alpha)(\alpha-\gamma)|^2+\delta^2|\tilde{N}(\alpha)|^2\big)^{\frac12},
		\end{aligned}
	\end{equation}
	for which we have the lower bound, also valid for $u_h$ \eqref{uh_tan},
		\begin{equation}\label{uhvh_lower_tan}
		\begin{aligned}
			|v_h|^2&\geq \frac12 |\pa Z|_{\text{inf}}^2\frac{|\pa Z|_{\text{inf}}^2}{\|\pa Z\|_{L^\infty}^2}\Big( |\alpha-\gamma|^2+\frac12\delta^2\Big).
		\end{aligned}
	\end{equation}
	Then, we perform the following splitting
	\begin{equation}\label{I5split_tan}
		\begin{aligned}
			I_5&=I_{5,1}+I_{5,2}+I_{5,3}+I_{5,4},
		\end{aligned}
	\end{equation}
	where
	\begin{equation*}
		\begin{aligned}
			I_{5,1}&=g(\alpha)h_\tau\cdot\pa Z_1(\alpha)\int_{B_\eta}\prod_{j=2}^3 ((\alpha-\gamma+h_\tau)\cdot\pa Z_j(\alpha)+\delta\tilde{N}_j(\alpha))\big(u_h^{-5}-v_h^{-5}\big)d\gamma,
		\end{aligned}
	\end{equation*}
	\begin{equation*}
		\begin{aligned}
			I_{5,2}&\!=\!g(\alpha)h_\tau\!\cdot\!\pa Z_2(\alpha)\!\!\int_{B_\eta}\!\!\!\!\!\!((\alpha\!-\!\gamma)\!\cdot\!\pa Z_1(\alpha)\!+\!\delta\tilde{N}_1(\alpha))((\alpha\!-\!\gamma\!+\!h_\tau\!)\cdot\!\pa Z_3(\alpha)\!+\!\delta\tilde{N}_3(\alpha))\big(u_h^{-5}\!-\!v_h^{-5}\big)d\gamma,
		\end{aligned}
	\end{equation*}
	\begin{equation*}
		\begin{aligned}
			I_{5,3}&=g(\alpha)h_\tau\cdot\pa Z_3(\alpha)\int_{B_\eta}\prod_{j=1}^2 ((\alpha-\gamma)\cdot\pa Z_j(\alpha)+\delta\tilde{N}_j(\alpha))\big(u_h^{-5}-v_h^{-5}\big)d\gamma,
		\end{aligned}
	\end{equation*}
	\begin{equation*}
		\begin{aligned}
			I_{5,4}&=g(\alpha)\int_{B_\eta}\prod_{j=1}^3 ((\alpha-\gamma)\cdot\pa Z_j(\alpha)+\delta\tilde{N}_j(\alpha))\big(u_h^{-5}-v_h^{-5}-u_0^{-5}+v_0^{-5}\big)d\gamma.
		\end{aligned}
	\end{equation*}
	We bound $I_{5,1}$ as follows
	\begin{equation*}
		\begin{aligned}
			|I_{5,1}|&\leq \|g\|_{L^\infty}\|\pa Z\|_{L^\infty}^3|h_\tau|\int_{B_\eta}(|\alpha-\gamma+h_\tau|+\delta)^2 |u_h^{-5}-v_h^{-5}|d\gamma.
		\end{aligned}
	\end{equation*}
	We have that
	\begin{equation}\label{uh_vh_tan}
		\begin{aligned}
			u_h^2-v_h^2&=|Z(\alpha)-Z(\gamma)-(\alpha-\gamma)\cdot\pa Z(\alpha)|^2\\
			&\quad+2\big((\alpha-\gamma+h_\tau)\cdot\pa Z(\alpha)+\delta\tilde{N}(\alpha)\big)\cdot\big(Z(\alpha)-Z(\gamma)-(\alpha-\gamma)\cdot\pa Z(\alpha)\big)\\
			&\leq C\big(|\alpha-\gamma|+|\alpha-\gamma+h_\tau|+\delta\big)\|\pa Z\|_{L^\infty}\|\pa Z\|_{\dot{C}^{\sigma}}|\alpha-\gamma|^{1+\sigma},
		\end{aligned}
	\end{equation}
	thus
	\begin{equation*}
		\begin{aligned}
			|u_h^{-5}-v_h^{-5}|&\leq C\frac{\|\pa Z\|_{L^\infty}^7}{|\pa Z|_{\inf}^{14}} \frac{|u_h^2-v_h^2|}{\big(|\alpha-\gamma|^2+\frac12\delta^2\big)^{\frac72}}\\
			&\leq C \frac{\|\pa Z\|_{L^\infty}^8\|\pa Z\|_{\dot{C}^{\sigma}}}{|\pa Z|_{\inf}^{14}}\frac{\big(|\alpha-\gamma|+|\alpha-\gamma+h_\tau|+\delta\big)|\alpha-\gamma|^{1+\sigma}}{\big(|\alpha-\gamma|^2+\frac12\delta^2\big)^{\frac72}}.
		\end{aligned}
	\end{equation*}
	Introducing this bound back, we obtain that
	\begin{equation}\label{I51bound_tan}
		\begin{aligned}
			|I_{5,1}|&\leq C\frac{\|g\|_{L^\infty}\|\pa Z\|_{L^\infty}^{11}\|\pa Z\|_{\dot{C}^{\sigma}}}{|\pa Z|_{\inf}^{14}}|h_\tau|\int_{B_\eta}\frac{\big(|\alpha-\gamma|+|\alpha-\gamma+h_\tau|+\delta\big)^3|\alpha-\gamma|^{1+\sigma}}{\big(|\alpha-\gamma|^2+\frac12\delta^2\big)^{\frac72}} d\gamma\\
			&\leq C\frac{\|g\|_{L^\infty}\|\pa Z\|_{L^\infty}^{10+\sigma}}{|\pa Z|_{\inf}^{13+\sigma}}\|\pa Z\|_{\dot{C}^{\sigma}}|h_\tau|^\sigma.
		\end{aligned}
	\end{equation}
	The same bound holds for the terms $I_{5,2}$ and $I_{5,3}$.
	We are left with $I_{5,4}$.
	We split it further as follows
	\begin{equation}\label{I54split_tan}
		\begin{aligned}
			I_{5,4}=J_1+J_2,
		\end{aligned}
	\end{equation}
	with
	\begin{equation*}
		\begin{aligned}
			J_1&=g(\alpha)\int_{B_\eta}\prod_{j=1}^3 ((\alpha-\gamma)\cdot\pa Z_j(\alpha)+\delta\tilde{N}_j(\alpha))G(u_h,u_0)\big(u_0^2-u_h^2-(v_0^2-v_h^2)\big)d\gamma,
		\end{aligned}
	\end{equation*}
	\begin{equation*}
		\begin{aligned}
			J_2&=g(\alpha)\int_{B_\eta}\prod_{j=1}^3 ((\alpha-\gamma)\cdot\pa Z_j(\alpha)+\delta\tilde{N}_j(\alpha))\big(G(u_h,u_0)-G(v_h,v_0)\big)\big(v_0^2-v_h^2\big)d\gamma,
		\end{aligned}
	\end{equation*}
	where $G$ is defined in \eqref{Guhvh}.
	From \eqref{vh_tan} we have that
	\begin{equation}\label{vhv0_tan}
		\begin{aligned}
			v_0^2-v_h^2&=-|h_\tau\cdot\pa Z(\alpha)|^2-2(h_\tau\cdot\pa Z)\cdot(\pa Z(\alpha)(\alpha-\gamma)),
		\end{aligned}
	\end{equation}
	which, together with \eqref{uhu0_tan}, provides that
	\begin{equation*}
		\begin{aligned}
			J_1&=-2g(\alpha)\int_{B_\eta}\prod_{j=1}^3 ((\alpha-\gamma)\cdot\pa Z_j(\alpha)+\delta\tilde{N}_j(\alpha))G(u_h,u_0)\\
			&\quad\times\big(Z(\alpha)-Z(\gamma)-\pa Z(\alpha)(\alpha-\gamma)\big)\cdot\big(h_\tau \cdot\pa Z(\alpha)\big)d\gamma.
		\end{aligned}
	\end{equation*}
	Therefore, we obtain 
	\begin{equation}\label{J1bound_tan}
		\begin{aligned}
			|J_1|\leq \frac{\|g\|_{L^\infty}\|\pa Z\|_{L^\infty}^{10+\sigma}}{|\pa Z|_{\inf}^{13+\sigma}}\|\pa  Z\|_{\dot{C}^\sigma}|h_\tau|^\sigma.
		\end{aligned}
	\end{equation}
	We proceed to estimate $J_2$. We further split this term
	\begin{equation}\label{J2split_tan}
		J_{2}=\sum_{k=1}^6J_{2,k}    
	\end{equation}
	where 
	\begin{equation*}
		\begin{aligned}
			J_{2,k}&=g(\alpha)\int_{B_\eta}\pa Z_1(\alpha)\cdot(\alpha-\gamma)\pa Z_2(\alpha)\cdot(\alpha-\gamma)\pa Z_3(\alpha)\cdot(\alpha-\gamma)\\
			&\hspace{2cm}\times \big(\frac1{u_h^{6-k}u_0^k}-\frac1{v_h^{6-k}v_0^k}\big)\frac{v_0^2-v_h^2}{u_h+u_0}d\gamma,
		\end{aligned}
	\end{equation*}
	for $1\leq k\leq 5$, and
	\begin{equation*}
		\begin{aligned}
			J_{2,6}&=g(\alpha)\int_{B_\eta}\pa Z_1(\alpha)\cdot(\alpha-\gamma)\pa Z_2(\alpha)\cdot(\alpha-\gamma)\pa Z_3(\alpha)\cdot(\alpha-\gamma)(v_0^2-v_h^2)\\
			&\hspace{2cm}\times \big(\frac1{v_h^5v_0}+\frac1{v_h^4v_0^2}+\frac1{v_h^3v_0^3}+\frac1{v_h^2v_0^4}+\frac1{v_hv_0^5}\big)\big(\frac{1}{u_h+u_0}-\frac1{v_h+v_0}\big)d\gamma.
		\end{aligned}
	\end{equation*}
	To control $J_{2,1}$ a further splitting is given:
	$$
	J_{2,1}=\sum_{l=1}^6K_l
	$$
	where 
	\begin{equation*}
		\begin{aligned}
			K_{1}&=g(\alpha)\int_{B_\eta}\pa Z_1(\alpha)\cdot(\alpha-\gamma)\pa Z_2(\alpha)\cdot(\alpha-\gamma)\pa Z_3(\alpha)\cdot(\alpha-\gamma)\\
			&\hspace{2cm}\times \frac1{u_h^5} \big(\frac1{u_0}-\frac1{v_0}\big)\frac{v_0^2-v_h^2}{u_h+u_0}d\gamma,
		\end{aligned}
	\end{equation*}
	and
	\begin{equation*}
		\begin{aligned}
			K_{l}&=g(\alpha)\int_{B_\eta}\pa Z_1(\alpha)\cdot(\alpha-\gamma)\pa Z_2(\alpha)\cdot(\alpha-\gamma)\pa Z_3(\alpha)\cdot(\alpha-\gamma)\\
			&\hspace{2cm}\times \frac1{v_0u_h^{6-l}v_h^{l-2}}\big(\frac1{u_h}-\frac1{v_h}\big)\frac{v_0^2-v_h^2}{u_h+u_0}d\gamma
		\end{aligned}.
	\end{equation*}
	for $2\leq l\leq 6$. Estimate \eqref{uh_vh_tan} and \eqref{vhv0_tan} allows us to get
	\begin{equation*}
		\begin{aligned}
			|K_1|&\leq C\|g\|_{L^\infty}\frac{\|\partial_\alpha Z\|_{L^\infty}^{13}}{\|\partial_\alpha Z\|^{16}_{\text{inf}}}\|\partial_\alpha Z\|_{\dot{C}^{\sigma}}|h_\tau|\int_{B_\eta} \frac{|\alpha-\gamma|^{4+\sigma}(|\alpha-\gamma|+\delta)(|h_\tau|+|\alpha-\gamma|)}{(|\alpha-\gamma|^2+\frac12\delta^2)^{\frac92}} d\gamma\\
			&\leq  C \minspace \|g\|_{L^\infty}\frac{\|\pa Z\|_{L^\infty}^{12+\sigma}}{|\pa Z|_{\text{inf}}^{15+\sigma}}\|\pa Z\|_{\dot{C}^\sigma}|h_\tau|^\sigma.
		\end{aligned}
	\end{equation*}
	The rest of $K_l$ are bounded similarly. Hence $J_{2,1}$ is controlled. 
	The remaining terms  $J_{2,k}$, $k=2,...,6$, in \eqref{J2split_tan} are controlled in a similar manner to $J_{2,1}$, and all of them are bounded with the same bound (we omit the details to avoid repetition, as the estimates follow along the lines below \eqref{v0-u0}). Therefore, 
	\begin{equation}\label{J2bound_tan}
		\begin{aligned}
			|J_2|&\leq  C \minspace \|g\|_{L^\infty}\frac{\|\pa Z\|_{L^\infty}^{12+\sigma}}{|\pa Z|_{\text{inf}}^{14+\sigma}}\|\pa Z\|_{\dot{C}^\sigma}|h_\tau|^\sigma,
		\end{aligned}
	\end{equation}
	which together with \eqref{J1bound_tan} gives that
	\begin{equation}\label{I54bound_tan}
		|I_{5,4}|\leq  C \minspace \|g\|_{L^\infty}\frac{\|\pa Z\|_{L^\infty}^{12+\sigma}}{|\pa Z|_{\text{inf}}^{15+\sigma}}\|\pa Z\|_{\dot{C}^\sigma}|h_\tau|^\sigma.
	\end{equation}
	Joining this bound with the ones for $I_{5,1}$, $I_{5,2}$, $I_{5,3}$ \eqref{I51bound_tan} back in \eqref{I5split_tan} provides that
	\begin{equation}\label{I5bound_tan}
		|I_{5}|\leq  C \minspace \|g\|_{L^\infty}\frac{\|\pa Z\|_{L^\infty}^{12+\sigma}}{|\pa Z|_{\text{inf}}^{15+\sigma}}\|\pa Z\|_{\dot{C}^\sigma}|h_\tau|^\sigma.
	\end{equation}
	Finally, the term $I_6$ \eqref{Isplit_tangent} can be written as follows
	\begin{equation*}
		\begin{aligned}
			I_6&\!=\! g(\alpha)\!\int_{|w-\frac{h_\tau}{\delta}|\leq \frac{\eta}{\delta}} \frac{\prod_{j=1}^3(w\cdot \pa Z_j(\alpha)\!+\!\tilde{N}_j(\alpha))}{\big(|\pa Z(\alpha)w|^2\!+\!|\tilde{N}(\alpha)|^2\big)^{\frac52}}d\gamma-g(\alpha)\int_{|w|\leq \frac{\eta}{\delta}} \frac{\prod_{j=1}^3(w\cdot \pa Z_j(\alpha)\!+\!\tilde{N}_j(\alpha))}{\big(|\pa Z(\alpha)w|^2\!+\!|\tilde{N}(\alpha)|^2\big)^{\frac52}}d\gamma,
		\end{aligned}
	\end{equation*}
	and therefore 
	\begin{equation*}
		\begin{aligned}
			|I_6|&\leq C\|g\|_{L^\infty}\int_{-\pi}^\pi\Big(\int_{\frac{\eta}{\delta}-\frac{|h_\tau|}{\delta}}^{\frac{\eta}{\delta}}+\int_{\frac{\eta}{\delta}}^{\frac{\eta}{\delta}+\frac{|h_\tau|}{\delta}}\Big) \frac{\prod_{j=1}^3((\hat{w}\cdot\pa Z_j(\alpha))r+\tilde{N}_j(\alpha))}{\big(|\pa Z(\alpha)\hat{w}|^2r^2\!+\!|\tilde{N}(\alpha)|^2\big)^{\frac52}}r\minspace dr\\
			&\leq C\|g\|_{L^\infty}\frac{\|\pa Z\|_{L^\infty}^8}{|\pa Z|_{\inf}^{10}}\Big(\int_{\frac{\eta}{\delta}-\frac{|h_\tau|}{\delta}}^{\frac{\eta}{\delta}}+\int_{\frac{\eta}{\delta}}^{\frac{\eta}{\delta}+\frac{|h_\tau|}{\delta}}\Big)\frac{r(r^3+1)}{(r^2+1)^{\frac52}}dr\\
			&\leq C\|g\|_{L^\infty}\frac{\|\pa Z\|_{L^\infty}^8}{|\pa Z|_{\inf}^{10}}\frac{|h_\tau|}{|\eta|}.
		\end{aligned}
	\end{equation*}
	Together with \eqref{I1bound_tan}, \eqref{I2I3bound_tan}, \eqref{I4bound_tan}, and \eqref{I5bound_tan} in \eqref{Isplit_tangent}, we conclude that
	\begin{equation}\label{Ibound_tan}
		\begin{aligned}
			|I|\leq  C \bigg(&\frac{\|\pa Z\|_{L^\infty}^{12+\sigma}}{|\pa Z|_{\text{inf}}^{14+\sigma}}\Big(\|g\|_{\dot{C}^\sigma}+\|g\|_{L^\infty}\frac{\|\pa Z\|_{\dot{C}^\sigma}}{|\pa Z|_{\inf}}\Big)+\frac{\|g\|_{L^\infty}\|\pa Z\|_{L^\infty}^8}{|\pa Z|_{\inf}^{10}}|h_\tau|^{1-\sigma}\bigg)|h_\tau|^\sigma,
		\end{aligned}
	\end{equation}
	and hence, substituting \eqref{gbound},
	\begin{equation*}
	    \begin{aligned}
	        |S(f)(x+h_\tau\cdot\pa& Z(\alpha))-S(f)(x)|\\
	        &\leq C(1+|\partial D|) P(\|F(Z)\|_{L^\infty}\!+\!\|\pa Z\|_{L^\infty})\Big(\|f\|_{\dot{C}^{\sigma}}+\|f\|_{L^\infty}\|\pa Z\|_{\dot{C}^{\sigma}}\Big)|h_\tau|^\sigma.
	    \end{aligned}
	\end{equation*}
	
	This concludes the proof of Case 2. Together with \eqref{normal_bound} (Case 1) and recalling the splitting \eqref{normaltangentSplit}, the estimate for two points near the boundary in \textit{nearly} normal direction, i.e., assuming \eqref{htau_menor_delta}, is done.

	\vspace{.3cm}

	\subsubsection{{\em{\textbf{Regularity in {\em{nearly}} tangential direction:}}}} Here we consider the case $|h_\tau|\geq \frac14\frac{|\pa Z|_{\inf}}{\|\pa Z\|_{L^\infty}}\delta$. 
	Recalling the expression \eqref{xh_normal} for $x+h$ in this case, we can write
	\begin{equation*}
		\begin{aligned}
			S(f)(x+h)-S(f)(x)&=S(f)(Z(\alpha+\lambda)+\mu\tilde{N}(\alpha+\lambda))-S(f)(Z(\alpha+\lambda))\\
			&\quad+S(f)(Z(\alpha+\lambda))-S(f)(Z(\alpha))\\
			&\quad+S(f)(Z(\alpha))-S(f)(Z(\alpha)+\delta\tilde{N}(\alpha)).
		\end{aligned}
	\end{equation*}
	Given the bound \eqref{lambdabound} and that we are in the case $\delta<4\frac{\|\pa Z\|_{L^\infty}}{|\pa Z|_{\inf}}|h_\tau|$, we can apply the previous H\"older estimates for $S(f)$ along the surface \eqref{bound_boundary} and on the normal direction \eqref{normal_bound} to obtain that
	\begin{multline}\label{tangentialestimate}
			|S(f)(x+h)-S(f)(x)|\leq \\
			C(1+|\partial D|) P(\|F(Z)\|_{L^\infty}\!+\!\|\pa Z\|_{L^\infty})\big(\|f\|_{C^{\sigma}}+\|f\|_{L^\infty}\|\pa Z\|_{\dot{C}^{\sigma}}\big)|h|^\sigma.
	\end{multline}

\vspace{0.3cm}
	
	\subsection{Regularity far from the boundary}\label{reg_far}
	Consider two points $x$ and $x+h$ in $ D$ (or analogously both in $\R^3\smallsetminus  D$) such that they are sufficiently far from the boundary. That is, recalling \eqref{delta_cutoff}, we consider now that
	\begin{equation*}
		\min\{d(x,\partial D), d(x+h,\partial D)\}\geq L,\qquad    |h|\leq \frac{L}{2},
	\end{equation*}
	where
	\begin{equation*}
		L=\frac{1}{6}\Big(\frac{|\pa Z|_{\inf}}{18\|\pa Z\|_{\dot{C}^\sigma}}\Big)^{\frac{1}{\sigma}}\Big(\frac{\|\pa Z\|_{L^\infty}}{|\pa Z|_{\inf}}\Big)^{\frac12}.
	\end{equation*}
	Then, for any $y\in\partial D$ and $s\in(0,1)$,
	\begin{equation*}
		|x-y+sh|\geq |x-y|-|h|\geq \frac{|x-y|}{2},
	\end{equation*}
	and thus the mean value theorem gives that
	\begin{equation*}
		\begin{aligned}
			|k(x+h-y)-k(x-y)|&\leq C\frac{|h|}{|x-y|^3},\\
			|k(x+h-y)-k(x-y)|&\leq \frac{C}{|x-y|^2}.
		\end{aligned}
	\end{equation*}
	Hence, it follows that
	\begin{equation*}
		\begin{aligned}
			|S(f)(x+h)-S(f)(x)|&=\Big|\int_{\partial  D}\big(k(x+h-y)-k(x-y)\big) f(y) dS(y)\Big|\\ &\leq C\|f\|_{L^\infty}\int_{\partial  D}\frac{|h|^\sigma}{|x-y|^{2+\sigma}}dS(y)\\
			&\leq C\|f\|_{L^\infty}\|\pa Z\|_{L^\infty}|h|^\sigma\int_L^\infty \frac{dr}{r^{1+\sigma}},
		\end{aligned}
	\end{equation*}
	and we conclude that
	\begin{equation*}
		\begin{aligned}
			|S(f)(x+h)-S(f)(x)|&\leq C\frac{\|f\|_{L^\infty}\|\pa Z\|_{L^\infty}}{L^\sigma}|h|^\sigma\\
			&\leq C\|f\|_{L^\infty}\frac{\|\pa Z\|_{L^\infty}}{|\pa Z|_{\inf}}\frac{|\pa Z|_{\inf}^{\frac{\sigma}{2}}}{\|\pa Z\|_{L^\infty}^{\frac{\sigma}{2}}}\|\pa Z\|_{\dot{C}^\sigma}|h|^\sigma.
		\end{aligned}
	\end{equation*}
	
	\qed

	\subsection*{Acknowledgements}
	The authors would like to gratefully thank the referee for their careful reading of the manuscript. 	FG and EGJ were partially supported by the ERC through the Starting Grant project H2020-EU.1.1.-639227. FG was partially supported by the grant EUR2020-112271 (Spain).
	EGJ was partially supported by the ERC Starting Grant ERC-StG-CAPA-852741. This project has received funding from the European Union’s Horizon 2020 research and innovation programme under the Marie Skłodowska-Curie grant agreement CAMINFLOW No 101031111.

	\bibliographystyle{plain}
	\bibliography{references9}

\begin{thebibliography}{10}

\bibitem{BGLV2016}
A.~Bertozzi, J.~Garnett, T.~Laurent, and J.~Verdera.
\newblock The regularity of the boundary of a multidimensional aggregation
  patch.
\newblock {\em SIAM J. Math. Anal.}, 48(6):3789--3819, 2016.

\bibitem{BertozziConstantin1993}
A.~L. Bertozzi and P.~Constantin.
\newblock Global regularity for vortex patches.
\newblock {\em Comm. Math. Phys.}, 152(1):19--28, 1993.

\bibitem{CMOV2021}
J.~C. Cantero, J.~Mateu, J.~Orobitg, and J.~Verdera.
\newblock Regularity of the boundary of vortex patches for some non-linear
  transport equations.
\newblock {\em To appear in Anal. PDE, arXiv:2103.05356}, 2021.

\bibitem{CCFGG19}
\'Angel Castro, Diego C\'{o}rdoba, Charles Fefferman, Francisco Gancedo, and
  Javier G\'{o}mez-Serrano.
\newblock Splash singularities for the free boundary {N}avier-{S}tokes
  equations.
\newblock {\em Ann. PDE}, 5(1):Paper No. 12, 117, 2019.

\bibitem{Chemin1993}
Jean-Yves Chemin.
\newblock Persistance de structures g\'{e}om\'{e}triques dans les fluides
  incompressibles bidimensionnels.
\newblock {\em Ann. Sci. \'{E}cole Norm. Sup. (4)}, 26(4):517--542, 1993.

\bibitem{CordobaGancedo2010}
Diego C\'{o}rdoba and Francisco Gancedo.
\newblock Absence of squirt singularities for the multi-phase {M}uskat problem.
\newblock {\em Comm. Math. Phys.}, 299(2):561--575, 2010.

\bibitem{CS19}
Daniel Coutand and Steve Shkoller.
\newblock On the splash singularity for the free-surface of a {N}avier-{S}tokes
  fluid.
\newblock {\em Ann. Inst. H. Poincar\'{e} Anal. Non Lin\'{e}aire},
  36(2):475--503, 2019.

\bibitem{Depauw1999}
Nicolas Depauw.
\newblock Poche de tourbillon pour {E}uler 2{D} dans un ouvert \`a bord.
\newblock {\em J. Math. Pures Appl. (9)}, 78(3):313--351, 1999.

\bibitem{Evans}
Lawrence~C. Evans.
\newblock {\em Partial differential equations}, volume~19 of {\em Graduate
  Studies in Mathematics}.
\newblock American Mathematical Society, Providence, RI, second edition, 2010.

\bibitem{GancedoG-J2018}
Francisco Gancedo and Eduardo Garc\'{i}a-Ju\'{a}rez.
\newblock Global regularity of 2{D} density patches for inhomogeneous
  {N}avier-{S}tokes.
\newblock {\em Arch. Ration. Mech. Anal.}, 229(1):339--360, 2018.

\bibitem{GancedoG-J2020}
Francisco Gancedo and Eduardo Garc\'{\i}a-Ju\'{a}rez.
\newblock Regularity results for viscous 3{D} {B}oussinesq temperature fronts.
\newblock {\em Comm. Math. Phys.}, 376(3):1705--1736, 2020.

\bibitem{GancedoG-J2021}
Francisco Gancedo and Eduardo Garc\'{i}a-Ju\'{a}rez.
\newblock {G}lobal regularity of {2D} {N}avier-{S}tokes free boundary with
  small viscosity contrast.
\newblock {\em To appear in Ann. Inst. H. Poincar\'e Anal. Non Lin\'eaire,
  arXiv:2109.08764}, 2021.

\bibitem{MajdaBertozzi2002}
Andrew~J. Majda and Andrea~L. Bertozzi.
\newblock {\em Vorticity and incompressible flow}, volume~27 of {\em Cambridge
  Texts in Applied Mathematics}.
\newblock Cambridge University Press, Cambridge, 2002.

\bibitem{MateuOrobitgVerdera2009}
Joan Mateu, Joan Orobitg, and Joan Verdera.
\newblock Extra cancellation of even {C}alder\'{o}n-{Z}ygmund operators and
  quasiconformal mappings.
\newblock {\em J. Math. Pures Appl. (9)}, 91(4):402--431, 2009.

\bibitem{MOV2011}
Joan Mateu, Joan Orobitg, and Joan Verdera.
\newblock Estimates for the maximal singular integral in terms of the singular
  integral: the case of even kernels.
\newblock {\em Ann. of Math. (2)}, 174(3):1429--1483, 2011.

\bibitem{Mitrea2Verdera2016}
Dorina Mitrea, Marius Mitrea, and Joan Verdera.
\newblock Characterizing regularity of domains via the {R}iesz transforms on
  their boundaries.
\newblock {\em Anal. PDE}, 9(4):955--1018, 2016.

\bibitem{PaicuZhang2020}
Marius Paicu and Ping Zhang.
\newblock Striated regularity of 2-{D} inhomogeneous incompressible
  {N}avier-{S}tokes system with variable viscosity.
\newblock {\em Comm. Math. Phys.}, 376(1):385--439, 2020.

\bibitem{Stein1970}
Elias~M. Stein.
\newblock {\em Singular integrals and differentiability properties of
  functions}.
\newblock Princeton Mathematical Series, No. 30. Princeton University Press,
  Princeton, N.J., 1970.

\bibitem{Stein93}
Elias~M. Stein.
\newblock {\em Harmonic analysis: real-variable methods, orthogonality, and
  oscillatory integrals}, volume~43 of {\em Princeton Mathematical Series}.
\newblock Princeton University Press, Princeton, NJ, 1993.

\bibitem{Torchinsky1986}
Alberto Torchinsky.
\newblock {\em Real-variable methods in harmonic analysis}, volume 123 of {\em
  Pure and Applied Mathematics}.
\newblock Academic Press, Inc., Orlando, FL, 1986.

\end{thebibliography}

\end{document}